\documentclass[12pt]{article}
\usepackage{amsmath}
\usepackage{amssymb}

\newcommand{\M}{{\cal M}}
\newcommand{\htilde}{{\tilde h}}
\newcommand{\Htilde}{{\tilde H}}
\newcommand{\supp}{{\rm supp}}
\newcommand{\algebra}{{\cal A}}
\newcommand{\alexr}{{\cal A}_R}
\newcommand{\anbv}{{\cal A}_{\nbv}}
\newcommand{\alexc}{{\cal A}_C}
\newcommand{\balexr}{{\cal B}_R}
\newcommand{\balexc}{{\cal B}_C}
\newcommand{\Rbar}{\overline{\R}}
\newcommand{\Partition}{{\cal P}}
\newcommand{\D}{{\cal D}}
\newcommand{\Dp}{{\cal D}'}
\newcommand{\nbv}{{\cal NBV}}
\newcommand{\nbvl}{{\cal NBV}_\lambda}

\newcommand{\ebv}{{\cal EBV}}
\newcommand{\bv}{{\cal BV}}
\newcommand{\intinf}{\int^\infty_{-\infty}}

\newcommand{\N}{{\mathbb N}}
\newcommand{\R}{{\mathbb R}}

\newcommand{\fn}{\!:\!}
\newcommand{\lsum}{\sum\limits}

\newcommand{\lsup}{\sup\limits}

\newcommand{\qed}{\mbox{$\quad\blacksquare$}}
\newcounter{examples}
\renewcommand{\theexamples}{(\alph{examples})\mbox{ }}
\newtheorem{theorem}{Theorem}

\newtheorem{prop}[theorem]{Proposition}
\newtheorem{corollary}[theorem]{Corollary}
\newtheorem{remark}[theorem]{Remark}
\newtheorem{defn}[theorem]{Definition}
\newtheorem{example}[theorem]{Example}
\newtheorem{proposition}[theorem]{Proposition}
\newtheorem{convention}[theorem]{Convention}
\begin{document}
\hspace{-2cm}
\raisebox{12ex}[1ex]{\fbox{{\footnotesize
Preprint
November 14, 2009.\quad To appear in 
{\it Illinois Journal of Mathematics.}
}}}

\begin{center}
{\large\bf The regulated primitive integral}
\vskip.25in
Erik Talvila\footnote{Supported by the
Natural Sciences and Engineering Research Council of Canada.
}\\ [2mm]
{\footnotesize
Department of Mathematics and Statistics \\
University of the Fraser Valley\\
Abbotsford, BC Canada V2S 7M8\\
Erik.Talvila@ufv.ca}
\end{center}

{\footnotesize
\noindent
{\bf Abstract.} 
A function on the real line  is called regulated if 
it has  a left limit and a  right limit at each point.
If $f$ is a Schwartz distribution on the real line
such that $f=F'$ (distributional or weak derivative) for
a regulated function $F$ then the regulated primitive
integral
of $f$ is $\int_{(a,b)}f=F(b-)-F(a+)$, with similar 
definitions for other types of intervals.
The space of
integrable distributions is a Banach space and Banach
lattice under the Alexiewicz norm.
It contains the
spaces of Lebesgue and Henstock--Kurzweil integrable
functions as  continuous embeddings.
It is the
completion
of the space of signed Radon measures in the Alexiewicz
norm.   Functions of
bounded variation form the dual space and the space of
multipliers.
The integrable
distributions are a module over the functions of bounded variation.  
Properties such
as integration by parts, change of variables, H\"older inequality,
Taylor's theorem and convergence
theorems are proved. 
\\
{\bf 2000 subject classification:} 26A39, 46E15, 46F05, 46G12
}\\

\section{Introduction}\label{introduction}
One way of defining an integral is via its primitive.  The primitive
is a function whose  derivative  is in some  sense equal
to the  integrand.  For example, if  $f$ and $F$ are  functions
on the real  line and  $F$ is absolutely continuous ($AC$) such that
$F'(x)=f(x)$
for almost all $x\in(a,b)$, then the  Lebesgue integral of $f$ is
  $\int_a^bf=F(b)-F(a)$.  The
same  method can be used to define the Henstock--Kurzweil integral,
for which $F\in ACG*$.  We  get  the wide Denjoy integral if
$F\in ACG$  and we use the approximate derivative. See, for example,
\cite{saks} for the definitions of these function spaces and the
wide Denjoy integral.  The Henstock--Kurzweil integral is equivalent
to the Denjoy integral  and is defined
further in this Introduction.  Note
that  $C^1\subsetneq AC\subsetneq ACG*\subsetneq ACG\subsetneq C^0$.
The symbols $\subset$ and $\supset$ allow set equality.
If we use the distributional derivative then the  primitives 
need not have any pointwise differentiation
properties.  See \cite{talviladenjoy} for an integral with
continuous functions as primitives.  In this paper we will describe
an integral  whose primitives are regulated  functions.  This integral will
contain all of the integrals above.  As well,  it  can integrate
signed Radon measures. 

We  will now  describe the space of primitives for the regulated  primitive
integral.  Function $F\fn\R\to\R$ is  {\it regulated} if it has left and right 
limits  at  each  point in $\R$.  For  $a\in\R$ write $F(a-)=\lim_{x\to a^-}
F(x)$, $F(a+)=\lim_{x\to a^+}F(x)$, $F(-\infty)=\lim_{x\to-\infty}F(x)$,
$F(\infty)=\lim_{x\to\infty}F(x)$.  Then  $F$  is left continuous  at
$a\in\R$ if $F(a)=F(a-)$ and right  continuous if  $F(a)=F(a+)$.  We
define $\balexr=\{F\fn\R\to\R\mid F \text{ is regulated and 
left continuous on } \R, 
F(-\infty)=0,
F(\infty)\in\R\}$. Hence, elements of $\balexr$ are real-valued functions
defined on the extended real line $\Rbar=[-\infty,\infty]$.
It  will  be shown  below that $\balexr$  is  a Banach
space under  the uniform norm.  The space of integrable  distributions
will  be those  distributions  that are  the  distributional derivative
of  a  function in  $\balexr$.  We will see that elementary properties of
distributions can be used to prove that the set of integrable
distributions is a Banach space isometrically isomorphic to  $\balexr$.
Most  of the usual properties of integrals hold: fundamental theorem of
calculus, addivity over intervals, integration by parts, change of
variables, H\"older inequality, Taylor's theorem,
convergence theorems.  The multipliers and  dual space  are the functions
of  bounded  variation.  This defines a product that makes the integrable
distributions into a module over the functions of bounded variation.
There is a  Banach lattice  structure.
The subspace of primitives of bounded variation corresponds to absolutely
convergent integrals.  Each integrable distribution is a finitely additive
measure defined on the algebra of sets of bounded variation.  We get a 
finite measure if and only if the convergence is absolute.  These 
distributions
are signed Radon measures.
The regulated functions are continuous in  the topology of half-open
intervals.
It is shown that the space of integrable distributions
is the completion of the space of signed Radon measures in the 
Alexiewicz norm.  See the paragraph preceding Theorem~\ref{alexrproperties}
for the definition.  This embedding is continuous.  Note, however, that
the  spaces of  Lebesgue and Henstock--Kurzweil integrable functions
are separable in  the Alexiewicz norm topology while our space of integrable 
distributions  is not
separable. Hence, it is not the  completion of these function spaces.

The integral in the present paper has several possible extensions
to Euclidean spaces.
In $\R^n$, geometrical considerations change the character of  the integral.
And there is the problem of which differential operator to invert.
There are integrals associated
with inverting the $n$th order distributional differential operator $\partial^n/\partial x_1\cdots
\partial x_n$.  For continuous primitives, this type of integral was introduced in  \cite{pmikusinski1}
and developed systematically in \cite{ang}.  At the other extreme,
there are integrals that invert the first order
 distributional divergence operator.  See \cite{pfeffer}.
If a set $S\subset\R^n$ has a normal vector at almost all points of its boundary
then we can use this
direction to define limits along the normal from within and without $S$.
This then defines functions that are regulated on the boundary of this
particular set $S$.  The divergence theorem for
sets of finite perimeter (\cite{evans}, \cite{ziemer}) can then be used to 
define an integral over $S$ 
for distributions that are the distributional derivative of a
regulated function.  If $S$ is a hypercube in $\R^n$ then its boundary is
a hypercube in $\R^{n-1}$ so that the divergence theorem in $\R^n$
yields the integral of \cite{ang}. Details will be published elsewhere.

If $\mu$ is a Borel measure on the real line we use the notation $L^p(\mu)$ 
($1\leq p<\infty$)
for the functions $f\fn\R\to\R$ such that
the Lebesgue integral $\int_\R|f|^p\,d\mu$ exists.  For Lebesgue measure $\lambda$
we write $L^p$.

To
proceed further we  will fix some notation for distributions.
The {\it test functions} are $\D=C^\infty_c(\R)$,  i.e., the smooth functions
with compact support.  The {\it support} of  a  function $\phi$ is the closure
of  the  set on  which $\phi$ does not  vanish.  Denote this as $\supp(\phi)$.
There is a notion of  continuity in $\D$.  If
$\{\phi_n\}\subset\D$ then $\phi_n\to\phi\in\D$ if there is a compact
set $K\subset\R$  such  that for  all  $n\in\N$, $\supp(\phi_n)\subset K$,
and for each integer $m\geq 0$, $\phi^{(m)}_n\to \phi^{(m)}$ uniformly
on $K$ as $n\to\infty$.  The distributions are the  continuous linear
functionals on $\D$,  denoted $\Dp$.  If $T\in\Dp$ then $T\fn\D\to\R$
and we write $\langle T,\phi\rangle\in\R$ for $\phi\in\D$.  If
$\phi_n\to\phi$ in $\D$ then $\langle T,\phi_n\rangle\to\langle T,\phi\rangle$
in $\R$.  And, for all $a_1, a_2\in\R$ and all $\phi,\psi\in\D$, $\langle T,
a_1\phi+a_2\psi\rangle =a_1\langle T,\phi\rangle+a_2\langle T,\psi\rangle$.  If
$f\in L^p_{{\rm loc}}$  for some $1\leq p\leq\infty$ then $\langle T_f,\phi\rangle=\int_{-\infty}^\infty
f\phi$ defines a distribution. The  differentiation  formula $\langle T',
\phi\rangle=-\langle T,\phi'\rangle$ ensures that all distributions have
derivatives of all  orders which are themselves  distributions.  This is
known as the {\it distributional derivative} or {\it weak derivative}.  The
formula follows by mimicking integration
by parts 
in the case  of  $T_f$ where $f\in C^1$.  We  will  usually denote distributional
derivatives by $F'$ and  pointwise derivatives by $F'(t)$.
For $T\in\Dp$ and $t\in\R$  the {\it translation} $\tau_t$ is defined
by $\langle\tau_tT, \phi\rangle=\langle T, \tau_{-t}\phi\rangle$
where $\tau_t\phi(x)=\phi(x-t)$ for
$\phi\in\D$.   Most of the results on distributions we
 use
can  be found in \cite{friedlander}.

L. Schwartz \cite{schwartz} 
defined the integral of distribution $T$ as $\langle T, 1\rangle$,
provided this exists.  This agrees with $\intinf f$ if $T=T_f$ for $f\in L^1$
(i.e., $f$ is integrable with respect to Lebesgue measure).
In Schwartz's definition the integral is then a linear functional on
the constant functions.  We will see that, as  a result   of  the
H\"older inequality (Theorem~\ref{holder}),
our integral can be viewed as a continuous linear functional 
on functions of bounded variation (Theorem~\ref{dual}).  
It then extends Schwartz's definition.
Other methods of integrating distributions have been considered by
J.~Mikusi\'nski, J. A.~Musielak and R.~Sikorski.  See the last paragraph
in Section~11 of \cite{talviladenjoy} for references.
As will be seen in Section~\ref{Sectionintegrationbyparts} below, the
integration by parts formula connects the regulated primitive integral
with a type of Stieltjes integral that has been studied by J.~Kurzweil,
\v{S}.~Schwabik and M.~Tvrd\'y \cite{tvrdy1989}.

Denjoy's original constructive approach to a nonabsolute integral that
integrated all pointwise and approximate derivatives was a type of
transfinite induction applied to sequences of Lebesgue integrals.  The
method of $ACG*$ functions and $ACG$ functions is due to Lusin.
See \cite{jeffery} and \cite{saks}.
In what follows we define an integral that integrates weak derivatives
of regulated functions.  It also has a definition  based on a
previously defined integral;  in this case from Riemann
integrals on compact intervals.

Another approach to nonabsolute integration is through Riemann sums.
The Henstock--Stieltjes
integral is defined as  follows (\cite[\S 7.1]{mcleod} where it is called
the gauge integral).
A {\it gauge} on $\Rbar$ is  a  mapping $\gamma$  from $\Rbar$ to  the open
intervals in $\Rbar$ (cf. Remark~\ref{compactification}) with the  property that  for each $x\in\Rbar$, $\gamma(x)$
is an open interval containing $x$.  Note that this requires
$\gamma(\pm\infty)=\Rbar$ or $\gamma(-\infty)=
[-\infty,  a)$ or $\gamma(\infty)=(b,\infty]$ for some $a,b\in\R$.  A
{\it tagged partition} is a finite set of  pairs of closed intervals
and tag points in the extended real line, $\Partition=\{([x_{n-1},x_n],z_n)\}_{n=1}^N$ for some $N\in\N$ 
such that $z_n\in[x_{n-1},x_n]$ for each $1\leq n\leq N$ and
$-\infty=x_0<x_1<x_2<\ldots<x_N=\infty$.  In addition, $z_0=-\infty$
and $z_N=\infty$.  Given a gauge, $\gamma$, the partition $\Partition$ is said
to be $\gamma${\it -fine} if $[x_{n-1},x_n]\subset\gamma(z_n)$ for
each $1\leq n\leq N$.  If  $F,g\fn\R\to\R$ then $F$  is integrable
with  respect to $g$ if there is $A\in\R$ such that for  all
$\epsilon>0$ there is  gauge $\gamma$ such  that  for each
$\gamma$-fine tagged  partition we have $|\sum_{n=1}^NF(z_n)\left[
g(x_n)-g(x_{n-1})\right]-A|<\epsilon$.  We will use the Henstock--Stieltjes 
integral
only for  regulated functions $F$ and $g$ so we can use limits to define
the values of these functions at $\pm\infty$.  If  $g(x)=x$ we have the
Henstock--Kurzweil integral, i.e., integration with respect to Lebesgue
measure.  In this case, we take $F(\pm\infty)=0$.
For integration over a compact interval, a  function is
Riemann integrable if and  only  if the gauge $\gamma$  can  be taken
to be constant, i.e., 
$\gamma(x)=(x-\delta,x+\delta)$ for some constant $\delta>0$.

If  function $F$ has  a pointwise derivative
at each point in $[a,b]$ then the derivative  is integrable in the
Henstock--Kurzweil sense and $\int_a^b F'(x)\,dx=F(b)-F(a)$.
In this  sense, the Henstock--Kurzweil integral inverts the pointwise
derivative  operator.  It is well  known that the  Riemann and Lebesgue
integrals do  not have this property.  For
details see \cite{mcleod}.  
There are functions for which this  fundamental theorem  of  calculus
formula holds and yet these  functions do not have a pointwise
derivative at
each point.  In this sense the Henstock--Kurzweil integral is not
the   inverse  of the  pointwise  derivative.  The $C$-integral of
B.~Bongiorno, L.~ Di Piazza  and D.~Priess is defined
using Riemann sums and  a modification of  the  gauge process  above.
A function has  a $C$-integral if and only if it
is everywhere the  pointwise derivative of its primitive.
See \cite{bongiorno}.  In
this sense the  $C$-integral  is the inverse of the pointwise
derivative.  The  integral  defined in the present paper inverts
the distributional  derivative but only for primitives that are
regulated functions. 
The
restriction to regulated  primitives is useful as  it leads to  a
Banach space of integrable distributions.

\section{The regulated primitive integral}\label{rpi}
Define
$\alexr =\{f\in\Dp\mid f=F' \mbox{ for some }
F\in \balexr\}$.  A distribution $f$ is integrable if it is the distributional
derivative of some primitive $F\in\balexr$, i.e., for all $\phi\in\D$ we have
$\langle f,\phi\rangle
=\langle F',\phi\rangle=-\langle F,\phi'\rangle=-\int_{-\infty}^\infty
F(t)\phi'(t)\,dt$.  Since $F$ is regulated and $\phi$ is smooth with compact
support, the last integral exists as a Riemann integral.  We will use the
following convention in labeling primitives of elements in $\alexr$.
\begin{convention}\label{convention}
When  $f,g,f_1,\htilde,$ etc. are in  $\alexr$ we
will denote their respective primitives in $\balexr$ by $F,G,F_1,\Htilde,$ etc.
\end{convention}
It will
be shown in Theorem~\ref{alexrproperties} below that primitives are unique
and that the spaces $\alexr$  and $\balexr$ are  isometrically 
isomorphic, the  integral constituting a  linear isometry.
If $f\in\alexr$ and $-\infty
<a<b<\infty$ then
\begin{eqnarray}
\int_{(a,b)}f & = & \int_{a+}^{b-}f=F(b-)-F(a+) =F(b)-F(a+)\label{intdefn1}\\
\int_{(a,b]}f & = & \int_{a+}^{b+}f=F(b+)-F(a+)\label{intdefn2}\\
\int_{[a,b)}f & = & \int_{a-}^{b-}f=F(b-)-F(a-)=F(b)-F(a)\label{intdefn3}\\
\int_{[a,b]}f & = & \int_{a-}^{b+}f=F(b+)-F(a-)
=F(b+)-F(a).\label{intdefn4}
\end{eqnarray}
If $F$ is continuous then  these four  integrals agree.  For $a=-\infty$
and $b=\infty$ we write these four integrals as $\intinf f=F(\infty)$.
We can also define $\int_{\{a\}}f=\int_{[a,a]}f=\int_{a-}^{a+}f=F(a+)-F(a-)$.

Elements of $\balexr$ are tempered distributions of order one, while
elements of $\alexr$ are tempered distributions of order two.
See \cite{friedlander} for the definitions.

\section{Examples}\label{examples}
\stepcounter{examples}
\theexamples
If $F\in AC$ and $F'(t)=f(t)$ for almost  all $t\in\R$  then for
$\phi\in\D$ we  can integrate by parts to get
\begin{eqnarray*}
\langle F',\phi\rangle & = & -\langle F,\phi'\rangle = 
-\intinf F(t)\phi'(t)\,dt\\
 & = & \intinf F'(t)\phi(t)\,dt -\left[F(t)\phi(t)\right]_{t=-\infty}^\infty\\
 & = & \intinf f(t)\phi(t)\,dt = \langle f,\phi\rangle.
\end{eqnarray*}
Each of the integrals above is a Lebesgue integral.
It then follows that if $F(-\infty)=0$ and $F(\infty)$ exists then
$f$ is  Lebesgue integrable and $L^1\subsetneq \alexr$.   
Similarly, the regulated primitive integral contains the Henstock--Kurzweil
integral and wide Denjoy integral.  See \cite[p.~33, 34]{celidze} for  the
integration  by parts formula  for these integrals.  

\stepcounter{examples}
\theexamples
Suppose $F\in\balexr$ is continuous but differentiable nowhere.  Then
$f$ defined  by $f=F'\in\alexr$ and 
$\int_I f=F(b)-F(a)$ for $I=(a,b), (a,b], [a,b), [a,b]$
whenever  $-\infty\leq a < b\leq \infty$.  Note that 
for $\phi\in\D$ we have $\langle f,\phi\rangle = \langle F',\phi\rangle=
-\langle F,\phi'\rangle = -\intinf F(t)\phi'(t)\,dt$.  This  last  integral
exists  in  the Riemann  sense.

\stepcounter{examples}
\theexamples
If $F\fn\R\to \R$ is  a  continuous function such  that $F'(x)=0$ for
almost  all $x\in\R$ then for all $[a,b]\subset\R$ 
the Lebesgue integral $\int_a^b F'(t)\,dt$ exists  and is  zero, while
$\int_a^b F'=F(b)-F(a)$.  An  example of such a  function $F$ is
the Cantor--Lebesgue function (devil's staircase).

\stepcounter{examples}
\theexamples
Let $\bv$ denote the  functions of  {\it bounded  variation}, i.e., functions
$F$ for which $VF:=\sup\sum|F(x_i)-F(y_i)|$ is bounded, where the  supremum
is taken over all disjoint intervals  $\{(x_i,y_i)\}$. Note that if
$F\in\bv$ then $F$ is regulated and 
$F(\pm\infty)$ exist.
Hence, $F'\in\alexr$.  Although $F'(t)$ exists for
almost all $t\in\R$, and the Lebesgue integral $\int_a^b F'(t)\,dt$
exists,  it need  not equal  $F(b)-F(a)$. If $F\in L^1_{loc}$ then
the {\it essential variation} of $F$ is  
${\rm ess\,var} F:=\sup\intinf F(t)\phi'(t)\,dt$ where
the supremum is taken over  all $\phi\in\D$ with $\|\phi\|_\infty\leq 1$.
And, ${\rm ess\,var} F =\inf VG$ such that $F=G$ almost everywhere.
The essential  variation  can also be computed  by restricting the
points $x_i,  y_i$ above to be points of approximate continuity of $F$.
Denote the functions with bounded essential variation as $\ebv$.
If $F\in\ebv$ then the distributional derivative of $F$ is
a signed Radon  measure, i.e., there is a signed Radon  measure $\mu$ such that
for all $\phi\in\D$ we have  $\langle F',\phi\rangle
=-\langle F,\phi'\rangle = \intinf \phi(t)\,d\mu(t)$.  
Radon measures are Borel measures  that are finite  on compact
sets, inner regular with respect to compact sets ($\mu(E)=\sup\mu(K)$
where $E$ is a Borel set and the supremum is taken over all compact sets
$K\subset E$) and outer regular with respect to open  sets
($\mu(E)=\inf\mu(G)$
where $E$ is a Borel set and the infimum is taken over all open sets
$G\supset E$). 
In $\R$ the Radon measures are the Borel measures that are finite
on compact sets.  See, for example, \cite[\S 26]{bauer}.
A signed  Radon  measure  is then the  difference of
two finite Radon measures.  If  $\mu$ is a signed Radon measure then 
$F$ defined  by $F(x)=\int_{(-\infty,x)}d\mu$ is a function of  bounded
variation. For, if $\{(x_i,y_i)\}$ are disjoint intervals then
$\sum|F(x_i)-F(y_i)|=\sum\left|\int_{[x_i,y_i)}d\mu\right|\leq
\sum\int_{[x_i,y_i)}d|\mu|\leq|\mu|(\R)<\infty$.  
The regularity of $\mu$ shows $F\in\balexr$.
Hence, each signed Radon measure
is in $\alexr$.  Since functions of bounded variation can
have a pointwise derivative that vanishes almost everywhere, we cannot
use a descriptive definition of the integral of a measure
using the pointwise derivative of its primitive.

If $\nu$ is a Radon measure and $f\in L^1(\nu)$
then 
the set function $\mu$ defined by $\mu(E)=\int_Ef\,d\nu$
is a signed Radon measure. Hence, $\mu\in\alexr$.  If
$\nu$ is absolutely continuous 
with respect to Lebesgue
measure ($\nu\ll\lambda$) and $f\in L^1(\nu)$ then
it follows  from  the Radon--Nikod\'ym theorem that
$fd\nu/d\lambda\in L^1\subset\alexr$.

\stepcounter{examples}
\theexamples
A distribution $T$ is said to be {\it positive} if $\langle T,\phi\rangle\geq 0$
whenever $\phi\in\D$ with $\phi\geq 0$.  It is known that positive
distributions correspond to Radon measures, i.e., 
$T\in\Dp$ is positive if and only if
there is a Radon measure $\mu$ such that for all $\phi\in\D$ we
have $\langle T,\phi\rangle =\intinf \phi(t)\,d\mu(t)$.  
For example, \cite[page~17]{vladimirov}. 
An example
of a positive distribution in $\balexr$ is the Dirac distribution.
Define the Heaviside step function by 
$H_1(x)=0$ for $x\leq 0$ and $H_1(x)=1$ for $x>0$.  The Dirac
distribution is then given by $\langle\delta,\phi\rangle =\phi(0)$ ($\phi
\in\D$). And, $\langle H_1',\phi\rangle =-\int_0^\infty\phi'(t)\,dt=\phi(0)$
so $H_1'=\delta$, $H_1\in\balexr$, $\delta\in\alexr$. We have
$\int_{(0,1)}\delta=\int_{(0,1)}H_1'=H_1(1-)-H_1(0+)=1-1=0$ while
$\int_{[0,1)}\delta=H_1(1-)-H_1(0-)=1-0=1$.  Define
$H_2(x)=0$ for $x< 0$ and $H_2(x)=1$ for $x\geq 0$. In $\D$, $H_1=H_2$
and  $H_1'=H_2'=\delta$.  Note that $H_2\not\in\balexr$
but $\int_I H_2'=\int_I H_1'$ for every interval $I\subset\R$.  This
discrepancy in $\balexr$ is discussed in Remark~\ref{equivalencerelation}
below.  Note also that $\delta$ is a Radon measure defined by
$\delta(E)=\chi_E(0)$ for all $E\subset\R$.  And, $\int_{\{0\}}\delta =1$.

\stepcounter{examples}
\theexamples
If $\{a_k\}$ is a sequence in $\R$ such that $\sum_{1}^\infty a_k$ converges
(absolutely or conditionally) then we can define a function $F\fn[0,\infty)\to
\R$ by $F(x)=\sum_{1}^n a_k$ if $x\in(n,n+1]$ for some $n\in\N$ and
$F(x)=0$ if $x\leq 1$.  Then $F$ is regulated, left continuous, $F(0)=0$
and $F(\infty)=\sum_{1}^\infty a_k$.
We have $F'=f$ where $f\in\alexr$ is the distribution
$f=\sum_{1}^\infty a_k(\tau_k\delta)$.  (See the
Introduction for the definition of translation.)  This gives
$\int_{[1,N]}f=\sum_{1}^N a_k$ for each $N\in\N$ and for $N=\infty$.  
Hence, integration in $\alexr$ includes series.

\stepcounter{examples}
\theexamples
Some finitely additive measures  are also in $\alexr$.  For example,
if $f(t)=\sin(t^2)$  define $F(x):=\int_{-\infty}^x f(t)\,dt$.
Then $F(I)=F(b)-F(a)$ for interval $I$ with endpoints $-\infty\leq a<
b\leq\infty$ defines a finitely additive measure on the algebra generated
by intervals on the real line.   And, $F\in\balexr$ with $F'=T_f$.  Since
the  integral converges conditionally, $F$ is  not countably additive.
For example, $\sum_{0}^\infty F([\sqrt{2n\pi},\sqrt{(2n+1)\pi}))=\infty$
while  $\sum_{1}^\infty F([\sqrt{(2n-1)\pi},\sqrt{2n\pi}))=-\infty$
but $F([0,\infty))=\int_0^\infty\sin(t^2)\,dt=\sqrt{\pi}/2^{3/2}$.
A similar  example is obtained with $f(t)=(d/dt)[t^2\sin(t^{-4})]$.

\stepcounter{examples}
\theexamples
If $F\fn\Rbar\to\R$ is any function then the Riemann--Stieltjes integral
$\int_a^bdF=F(b)-F(a)$ exists for  all $(a,b)\subset\R$.  The
Riemann--Stieltjes integral then contains the  regulated primitive integral.
We will see that $\alexr$ is  a useful restriction since it is a Banach space.
Below it will be shown that we can define $\intinf f\,dg$ under more
general conditions than can be done for the Riemann--Stieltjes or
Lebesgue--Stieltjes integrals.

\section{Properties of the integral}
First we have some properties of  our  space of primitives.

\begin{theorem}[Properties of $\balexr$]\label{theorembalexr}
{\rm (a)} If  $F\in\balexr$ then it is {\it uniformly
regulated}, i.e., for each  $\epsilon>0$ there exists $\delta>0$
such that for each $x\in \R$, if $y\in(x-\delta,x)$ then
$|F(x-)-F(y)|<\epsilon$ and if
$y\in(x, x+\delta)$ then  $|F(x+)-F(y)|<\epsilon$.  If
$y<1/\delta$ then $|F(y)|<\epsilon$.  And, if $y>1/\delta$ then
$|F(\infty)-F(y)|<\epsilon$.  Similarly, $F$ is {\it uniformly
left continuous}.
{\rm (b)} If $F\in\balexr$ then $F$ is bounded and has at most a
countable number  of discontinuities.
{\rm (c)} Using pointwise operations, $\balexr$ is a Banach space
under the uniform norm: $\|F\|_\infty
=\sup_{x\in\R}|F(x)|$, for $F\in\balexr$.
{\rm (d)} $\balexr$ is not separable.
\end{theorem}

\bigskip
\noindent
{\bf Proof:}
(a) Let $\epsilon>0$. There is $\alpha<0$ such that if $y\leq\alpha$
then $|F(y)|<\epsilon$.  For each $x\in\R$ there is $\eta_x>0$ such
that if $y\in(x-\eta_x,x]$ then $|F(y)-F(x)|<\epsilon$.  There is
$\gamma_x>0$ such if $y\in(x,x+\gamma_x)$ then $|F(y)-F(x+)|<\epsilon$.
There is $\beta>0$ such that if $y\geq \beta$ then $|F(y)-F(\infty)|<\epsilon$.
Let $\zeta_x=\min(\eta_x,\gamma_x)$.
The family of open intervals $\{(x-\zeta_x,x+\zeta_x)\}_{x\in\R}$ forms
an open cover of the compact interval $[\alpha,\beta]$.  There is then
a finite index set $J\subset\R$ such that $\{(x-\zeta_x,x+\zeta_x)\}_{x\in J}$ 
is again an open cover of $[\alpha,\beta]$.
Now let $\delta=\min(-1/\alpha,1/\beta, \min_{x\in J}\zeta_x)$.  Since 
$\delta>0$ this shows $F$ is uniformly regulated and uniformly left continuous.

(b) In (a) let $\epsilon=1$.  Then 
$$
|F(x)|\leq 1+\max(|F(\alpha)|,
\max_{x\in J}(|F(x)|,|F(x+)|),|F(\beta)|)<\infty.
$$
See \cite[p.~225]{mcleod} for a proof  that there are  
at most countably many points of discontinuity.

(c)
By (b), if $F\in\balexr$ then $F$ is bounded and measurable.
To prove $\balexr$  is a Banach space, first  note it is a linear
subspace of $L^\infty(\R)$ since $\balexr$ is  clearly closed under
linear combinations.  And,  if $F\in\balexr$ such that $\|F\|_\infty=0$
then $F(x)=0$ for  almost  all $x\in\R$.  But $F$ is left continuous so
if there were $b\in\R$ such that $F(b)>0$ then there  is an
interval $(a,b]$ on which $F$ is positive, which is a  contradiction,
so $F(x)=0$ for all $x\in\R$.  Positivity, homogeneity and  the triangle
inequality are  inherited  from $L^\infty(\R)$.  
To show $\balexr$  is complete,  
suppose  $\{F_n\}$
is a Cauchy sequence in $\balexr$.  Then $\{F_n\}$ is a Cauchy sequence in
$L^\infty(\R)$ so there is $F\in L^\infty(\R)$ such  that $\|F-F_n\|_\infty
\to 0$.  To show $F$ is left continuous, suppose $a\in\R$.  For $x<a$  and
$n\in\N$,
\begin{eqnarray}
|F(a)-F(x)| & \leq & |F(a)-F_n(a)| +|F_n(a)-F_n(x)|+|F_n(x)-F(x)|\notag\\
 & \leq & 2\|F-F_n\|_\infty+|F_n(a)-F_n(x)|.\label{complete}
\end{eqnarray}
Given $\epsilon>0$,  fix  $n$ large enough so that $\|F-F_n\|_\infty<
\epsilon/3$.  Then let $x\to a-$.  Hence, $F$ is left continuous on $\R$.
Using $|F(a)|\leq \|F-F_n\|_\infty + |F_n(a)|$ we see that $F(-\infty)=0$.
We can see that $F$ has a right limit at $a\in\R$ by taking $x,y>a$
and letting $x,y\to a+$ in
$|F(x)-F(y)|\leq 2\|F-F_n\|_\infty + \|F_n(x)-F_n(y)|$.
And, $F(\infty)$ is  seen  to  exist by letting $x,y\to\infty$
in  this inequality.  Therefore, $F\in\balexr$ and
the space is complete.

(d) To see that $\balexr$ is not separable, consider the family
of translations $\{\tau_t H_1\mid t\in\R\}$.
The function $H_1$ is defined in Example~\ref{examples}(e).
Given $0<\epsilon<1/2$, for each $t\in\R$ a dense subset of $\balexr$
would have to contain a function $F_t$ with $|F_t|<\epsilon$
on $(-\infty,t]$ and $|F_t-1|<\epsilon$ on $(t,\infty)$.  Hence,
no such dense set can be countable.
\qed\\

Further properties of regulated functions can be found in 
\cite{frankova} and \cite{honig}.

\begin{remark}\label{compactification}
{\rm
Note that the construction in  (a) gives a compactification of $\R$.  
A topological
base for $\Rbar$ consists of the usual open intervals $(a,b)$ with
$-\infty\leq
a<b\leq\infty$, as well as $[-\infty,a)$ with $-\infty<a\leq\infty$,
and $(a,\infty]$ with $-\infty\leq a<\infty$.  This makes $\Rbar$ 
into a compact Hausdorff space.  A different topology is
introduced in Section~\ref{measure}, under which all functions
in $\balexr$ are continuous.
}
\end{remark}

We now  present some of the basic properties of  the  integral.
One of the main results is that $\alexr$ is  a Banach  space
under the Alexiewicz  norm.  For $f\in\alexr$ this is defined
as $\|f\|=\|F\|_\infty$ where, as usual, $F$ is  the unique
primitive  in $\balexr$ (Convention~\ref{convention}).  Linear combinations are defined by
$\langle a_1f_1+a_2f_2,\phi\rangle=\langle a_1F_1'+a_2F_2',\phi\rangle$
for $\phi\in\D$; $a_1, a_2\in\R$; $f_1,f_2\in\alexr$ with primitives
$F_1,F_2\in\balexr$.
\begin{theorem}[Basic properties of the integral]\label{alexrproperties}
{\rm (a)} The integral  is  unique.
{\rm (b)} Addivity  over intervals.  If $f\in\alexr$ then for all $-\infty\leq a
<b<c<\infty$ we have $\int_{(a,b]}f+\int_{(b,c]}f=\int_{(a,c]}f$.  There
are similar formulas for other intervals.
{\rm (c)} With the Alexiewicz norm, $\alexr$ is  a  Banach
space.  The integral  provides a  linear isometry and isomorphism
between $\alexr$
and $\balexr$. 
{\rm (d)} $\alexr$ is not separable.
{\rm (e)} Linearity.  If $f_1, f_2\in\alexr$ and $a_1, a_2\in\R$
then $a_1f_1+a_2f_2\in\alexr$ and $\int_{-\infty}^\infty(a_1f_1+a_2f_2)
=a_1\int_{-\infty}^\infty f_1+ a_2\int_{-\infty}^\infty f_2$.
{\rm (f)} Reverse limits of integration.  Let $-\infty\leq a_1<a_2\leq \infty$
and $\epsilon_1, \epsilon_2\in\{+,-\}$.  Then
$\int_{a_1\epsilon_1}^{a_2\epsilon_2}f=-\int_{a_2\epsilon_2}^{a_1\epsilon_1}f$.
If $a_1=-\infty$ then we don't need $\epsilon_1$ and if
$a_2=\infty$ then we don't need $\epsilon_2$.
\end{theorem}

\bigskip
\noindent
{\bf Proof:}
{\rm (a)}
To prove the  integral is unique we need to  prove
primitives  in $\balexr$ are unique.   Suppose $F,G\in\balexr$ and
$F'=G'$.  Then $(F-G)'=0$  and the only solutions of  this  distributional
differential equation are the constant distributions \cite[\S2.4]{friedlander}.
The only constant distribution in $\balexr$ is  the zero function.

{\rm (b)} Note that 
$\left[F(b+)-F(a+)\right] + \left[F(c+)-F(b+)\right]=F(c+)-F(a+)$.

(c) Linearity of the distributional derivative shows $\alexr$ is a 
linear subspace of $\Dp$.  To prove $\|\cdot\|$ is a norm, let $f,g\in
\alexr$.

(i) By uniqueness of the primitive, $\|0\|=\|0\|_\infty=0$.  If $\|f\|=0$
then $\|F\|_\infty=\sup_{x\in\R}|F(x)|=0$ so $F(x)=0$ for all $x\in\R$ and therefore $f=F'=0$.

(ii) Let $k\in\R$.  Then $(kF)'=kF'$ so $\|kf\|=\|kF\|_\infty=|k|\|F\|_\infty
=|k|\,\|f\|$.

(iii) Since $f+g=F'+G'=(F+G)'$ we have $\|f+g\|=\|F+G\|_\infty\leq
\|F\|_\infty+\|G\|_\infty=\|f\|+\|g\|$.

This shows $\alexr$ is a normed linear space.  To prove it is complete,
suppose $\{f_n\}$ is  a Cauchy sequence in $\alexr$.  Then $\|F_n-F_m\|_\infty
=\|f_n-f_m\|$ so $\{F_n\}$ is a Cauchy  sequence
in $\balexr$.  There is $F\in\balexr$
such that $\|F_n-F\|_\infty\to 0$.  And then $\|f_n-F'\|=\|F_n-F\|_\infty
\to 0$.  Since $F\in\balexr$ we have $F'\in\alexr$
and $\alexr$ is complete.

A linear bijection $\psi\fn\alexr\to\balexr$ is given by $\psi(f)=F$
where $f\in\alexr$ and $F$ is its unique primitive in $\balexr$.
Since the integral is  linear, so is $\psi$.  It is an isometry
because $\|f\|=\|F\|_\infty=\|\psi(f)\|_\infty$.

(d) To show that $\alexr$ is not separable, consider the set
$\{\tau_t \delta\mid t\in\R\}$. 
Now proceed as in the proof of Theorem~\ref{theorembalexr}(d).
 
(e)
Since $a_1f_1+a_2f_2=(a_1F_1+a_2F_2)'$ we  have
$\int_{-\infty}^\infty(a_1f_1+a_2f_2)=(a_1F_1+a_2F_2)(\infty)
=a_1F_1(\infty)+a_2F_2(\infty)=a_1\int_{-\infty}^\infty f_1
+a_2\int_{-\infty}^\infty f_2$.

(f) $\int_{a_1\epsilon_1}^{a_2\epsilon_2}f=
F(a_2\epsilon_2)-F(a_1\epsilon_1) =-[F(a_1\epsilon_1)-F(a_2\epsilon_2)]
=-\int_{a_2\epsilon_2}^{a_1\epsilon_1}f$.
\qed\\

No space of integrable functions or distributions for which primitives
are continuous can be dense in $\alexr$.  If $G$ is a continuous
primitive then $\|G'-H_1'\|=\|G-H_1\|_\infty\geq 1/2$.  Thus, 
$L^1$ is not dense in $\alexr$, nor are the spaces of Henstock--Kurzweil
or wide Denjoy integrable functions.
The completion of these spaces in the Alexiewicz norm is
the Banach space
$\alexc=\{f\in\Dp\mid f=F' \text{ for some } F\in\balexc\}$,
where $\balexc=\{F\in C^0(\Rbar)\mid F(-\infty)=0 \text{ and }
F(\infty)\in\R\}$. If $f\in\alexc$ and $F'=f$ where $F$ is its
unique primitive in $\balexc$ then   the {\it continuous primitive 
integral} of $f$ is $\int_a^b f=F(b)-F(a)$.  This  integral is discussed
in \cite{talviladenjoy}, where further references can also be found.
Note that the spaces of Henstock--Kurzweil
and wide Denjoy integrable functions are barrelled but not complete
under the Alexiewicz norm.

\begin{remark}\label{equivalencerelation}
{\rm In defining $\balexr$ we have chosen the primitives to  be left
continuous.  This is convenient but somewhat arbitrary.  
If two regulated functions have  the  same left  and
right  limit  at  each point  then the functions can be different on
a countable set but will still  define the same distribution and thus
have the same distributional derivative.  This does not affect the
integral since it only depends on limits at endpoints of an interval
and  not on the value of  the  primitive at the endpoints.  
It is clear that an equivalence relation between such primitives
could be established, namely, $F\equiv G$ if and only if
$F(x-)=G(x-)$ for all $-\infty<x\leq\infty$ and  $F(x+)=G(x+)$
for all $-\infty\leq x<\infty$.
An advantage of using left
continuous functions rather than just regulated functions 
is that the norm on $\balexr$ can be taken as
$\|F\|_\infty=\sup_{x\in\R}|F(x)|$ rather than essential supremum.
This choice also affects the lattice operations in Section~\ref{lattice}.
Other obvious conventions are to take primitives that are right 
continuous or for which $F(x) =[F(x-)+F(x+)]/2$.
As pointed out in \cite{kaltenborn}, any normalising condition
$F(x)=(1-\lambda) F(x-)+\lambda F(x+)$ suffices for fixed $0\leq\lambda\leq 1$.
In Lebesgue and Henstock--Kurzweil integration we have equivalence classes
of functions that agree almost everywhere.  In $\alexr$ there are no
such equivalence classes, for two distributions are equal if they
agree on all test functions.  For example, if $f,g\in L^1_{loc}$ 
and $f=g$ almost everywhere then $T_f=T_g$.
}\qed
\end{remark}

Our definition of the integral builds in half of the fundamental theorem
of calculus.  The other half follows easily from uniqueness.

\begin{theorem}[Fundamental theorem of calculus]\label{fundamental}
{\rm (a)} Let $f\in\alexr$.  Define $G_1(x)=\int_{(-\infty,x)}f$.  Then
$G_1=F$ on $\R$ and $G_1'=f$ in $\D$.  
Define $G_2(x)=\int_{(-\infty,x]}f$.  Then $G_2$ is right
continuous, $G_2(-\infty)=0$, $G_2(\infty)$ exists and $G_2'=f$.
{\rm (b)} Let $G$ be a regulated function with limits at $\pm\infty$.
Then $G'\in\alexr$ and, for all  $x\in\R$, $\int_{(-\infty,x)}G'=G(x-)-G(-\infty)$ and
$\int_{(-\infty,x]}G'=G(x+)-G(-\infty)$.
\end{theorem}
\bigskip
\noindent
{\bf Proof:} 
(a) Since $f\in\alexr$ there is a unique function 
$F\in\balexr$ such that $F'=f$ and $G_1(x)=
\int_{(-\infty,x)}f=F(x-)=F(x)$
for all $x\in\R$.
For $x\in\R$ we have $G_2(x)=\int_{(-\infty,x]}f=F(x+)$.  It
follows that $G_2$ is right continuous.  And, $\lim_{x\to-\infty}G_2(x)
=\lim_{x\to-\infty}F(x+)=F(-\infty)=0$.  As well, $\lim_{x\to\infty}G_2(x)
=\lim_{x\to\infty}F(x+)=F(\infty)$.  Therefore, $G_2=F$ except perhaps
on a countable set.  They then define the same distribution and $G_2'=F'=f$.

(b) Define $F_1(x)=G(x-)-G(-\infty)$.  Then $F_1\in\balexr$
and $F_1'=G'$ so $G'\in\alexr$. Since $G((x-)-)=G(x-)$ 
we have $\int_{(-\infty,x)}G'=
F_1(x-)=G(x-)-G(-\infty)$.  As well, $G((x-)+)=G(x+)$ so $\int_{(-\infty,x]}G'=
F_1(x+)=G(x+)-G(-\infty)$.\qed\\
 
As with  the  Henstock--Kurzweil integral there are  no improper
integrals.

\begin{theorem}[Hake theorem]\label{haketheorem}
Suppose $f\in\Dp$ and $f=F'$ for some regulated function $F$.  If $F(-\infty)$
and $F(\infty)$ exist in $\R$ then $f\in\alexr$ and
$\intinf f=\lim_{x\to\infty}\int_{(0,x)} f +\lim_{x\to-\infty}\int_{(x,0]} f
=F(\infty)-F(-\infty)$.
\end{theorem}
There are similar versions of this theorem on compact intervals.

If $T$ is a distribution and $G\fn\R\to\R$ is an increasing $C^\infty$
bijection then for test function $\phi$ define $\psi\in\D$ by
$\psi=(\phi\circ G^{-1})/(G'\circ G^{-1})$.  The composition $T\circ G
\in\Dp$
is then defined by $\langle T\circ G,\phi\rangle = \langle T,\psi\rangle$.
This follows from the change of variables formula for integration of 
smooth functions.  See \cite[\S7.1]{friedlander}.

For the case at hand we can reduce the requirements on $G$.
If $F\in\balexr$ and the real line can be partitioned into a finite
number of intervals, on each of which  $G$ is monotonic (i.e., $G$ is
piecewise monotonic) then
$F\circ G$ is regulated so we can integrate its derivative.
This gives a change of variables formula.

\begin{theorem}[Change of variables]\label{theoremchangevariables}
{\rm (a)} Let $F\in\balexr$. For each point in $\R$ let $G\fn\R\to\R$ have
left and right limits with values in $\Rbar$.  Let $G$ be  piecewise
monotonic with $\lim_{\pm\infty}G$ existing in $\Rbar$.
Then $F\circ G$ is regulated on  $\R$ with
real limits at $\pm\infty$.  Define $(F'\circ G)G':=(F\circ G)'$.  \\

\noindent
{\rm (b)} Let $f\in\alexr$.  In addition to {\rm (a)}, assume $G$ is 
increasing,
left continuous, $\lim_{-\infty}G=-\infty$ and
$\lim_{\infty}G\in(-\infty,\infty]$.  Then $F\circ G\in\balexr$ and 
$\intinf(f\circ G)G'=\intinf (F\circ G)'= \int_{(-\infty,G(\infty))}f
=F(G(\infty)-)$.\\

\noindent
{\rm (c)}   Let $f\in\alexr$. Assume $G$ as in {\rm (a)}. 
Let $-\infty< a_1<a_2<\infty$.  For each $i\in\{1,2\}$, let
$\sigma_i,\epsilon_i\in\{+,-\}$.  Then
\begin{eqnarray}
\int_{a_1\epsilon_1}^{a_2\epsilon_2}(f\circ G)G' & = & 
\int_{a_1\epsilon_1}^{a_2\epsilon_2}(F\circ G)' =
\int_{G(a_1\epsilon_1)\sigma_1}^{G(a_2\epsilon_2)\sigma_2}f\notag\\
 & = & (F\circ G)(a_2\epsilon_2)-(F\circ G)(a_1\epsilon_1)\notag\\
 & = & 
 F(G(a_2\epsilon_2)\sigma_2)-F(G(a_1\epsilon_1)\sigma_1).\label{changeofvariables}
\end{eqnarray}
For each $i\in\{1,2\}$, $\sigma_i=\epsilon_i$ if $G$ is increasing on
an interval with endpoints $a_i$ and $a_i\epsilon_i\delta$ for
some $\delta>0$, and $\sigma_i\not=\epsilon_i$ if $G$ is decreasing on
an interval with endpoints $a_i$ and $a_i\epsilon_i\delta$ for
some $\delta>0$.  If $G(a_i\epsilon_i)=\pm\infty$ then we don't need 
$\sigma_i$.

If $a_1=-\infty$ then replace $a_1\epsilon_1$ with $-\infty$ in 
\eqref{changeofvariables}.
If $a_2=\infty$ then replace $a_2\epsilon_2$ with $\infty$ in 
\eqref{changeofvariables}.  If $G$ is increasing in a neighbourhood
of $-\infty$ then $\sigma_1=+$.
If $G$ is decreasing in a neighbourhood
of $-\infty$ then $\sigma_1=-$.
If $G$ is increasing in a neighbourhood
of $\infty$ then $\sigma_2=-$.
If $G$ is decreasing in a neighbourhood
of $\infty$ then $\sigma_2=+$.  If $G(-\infty)\in\{-\infty,\infty\}$ then
we don't need $\sigma_1$ and if $G(\infty)\in\{-\infty,\infty\}$ then
we don't need $\sigma_2$.\\

\noindent
{\rm (d)} Let $f\in\alexc$.
(See the paragraph preceding Remark~\ref{equivalencerelation} for the
definition.)  Let $G$ be regulated with $\lim_{\pm\infty}G$ existing
in $\Rbar$.  Then $F\circ G$ is regulated and
\begin{eqnarray*}
\int_{a_1\epsilon_1}^{a_2\epsilon_2}(f\circ G)G' & = & 
\int_{a_1\epsilon_1}^{a_2\epsilon_2}(F\circ G)' = 
\int_{G(a_1\epsilon_1)}^{G(a_2\epsilon_2)}f\notag\\
 & = & (F\circ G)(a_2\epsilon_2)-(F\circ G)(a_1\epsilon_1).
\end{eqnarray*}
The last integral exists as a continuous primitive integral
\cite{talviladenjoy}.
\end{theorem}

\bigskip
\noindent
{\bf Proof:}  (a)  Let $x\in\R$.  For small enough $\delta>0$, $G$ is
monotonic on intervals with endpoints $x$ and $x\pm\delta$.
Suppose $G$ is decreasing on $(x,x+\delta)$.
If $\lim_{y\to x^+}G(y)\in\R$ then for each $\nu>0$ there exists
$\eta(\nu)>0$ such that if $y\in(x,x+\eta(\nu))$ then
$G(y)\in(G(x+)-\eta(\nu),G(x+))$.  Since $\lim_{z\to G(x+)-}F(z)$ exists
we have that for every $\epsilon>0$ there is $\nu(\epsilon)>0$ such
that if $z\in(G(x+)-\nu(\epsilon),G(x+))$ then
$|F(z)-F(G(x+)-)|<\epsilon$.  To show $\lim_{z\to x+}(F\circ G)(z)$
exists, let $\epsilon>0$ and let 
$y\in(x,x+\eta(\nu(\epsilon)))$.  Then $|(F\circ G)(y)-F(G(x+)-)|<\epsilon$.
Other cases are similar with only minor modifications, including
showing left or right continuity of $F\circ G$.  Similarly in part (c).

(c) Suppose $-\infty<a<b<\infty$ and for some $\delta>0$ we have $G$ increasing
on $(a-\delta,a)$ and decreasing on $(b,b+\delta)$.  Then
$a_1=a$, $a_2=b$,
$\epsilon_1=\sigma_1=-$, $\epsilon_2=+$, $\sigma_2=-$ and we have
\begin{eqnarray*}
\int_{a-}^{b+}(f\circ G)G' & = & \int_{[a,b]}(F\circ G)' 
= (F\circ G)(b+)-(F\circ G)(a-)\\
 & = & F(G(b+)-)-F(G(a-)-)
= \int_{G(a-)-}^{G(b+)-}f.
\end{eqnarray*}
Other cases are similar and (b) is included in (c).

(d) There need be no interval on which $G$ is monotonic.  However,
since $F$ is continuous we have
\begin{eqnarray*}
\int_{a_1\epsilon_1}^{a_2\epsilon_2}(f\circ G)G' & = & 
(F\circ G)(a_2\epsilon_2) - (F\circ G)(a_1\epsilon_1)\\
 & = & \lim_{x\to G(a_2\epsilon_2)}F(x) -  \lim_{x\to G(a_1\epsilon_1)}F(x) \\
 & = & F(G(a_2\epsilon_2)) - F(G(a_1\epsilon_1))\\
 & = & \int_{G(a_1\epsilon_1)}^{G(a_2\epsilon_2)}f.
\end{eqnarray*}
The cases when $a_1=-\infty$ or $a_2=\infty$ are similar.\qed\\

Note that in (c) there are 16 cases, depending on whether
$G$ is increasing or decreasing at each of the endpoints for the
four types of integrals in \eqref{intdefn1}.  There are four cases
for endpoints at $\pm\infty$. 

Note that $G$ need not be strictly monotonic but then we have to
use the left continuity of $F$ to interpret the integral.  For
example, if $f\in\alexr$ and $G=H_1$ then
$\intinf(F\circ H_1)'=\int_{-\infty}^\infty (f\circ H_1)\delta=(F\circ H_1)(\infty)-
(F\circ H_1)(-\infty)=F(1)-F(0)=F(1-)-F(0-)=\int_{[0,1)}f$.
And,
$G$ need not be bounded.  For example, let $G(x)=1+x^{-2}$ for $x\not=0$.
The value of $G$ at $0$ is immaterial.  Let $f\in\alexr$.  Then
$\int_{(-\infty,0)}(f\circ G)G' =(F\circ G)(0-)-(F\circ G)(-\infty)
=F(\infty)-F(1+)$.  We have $a_1=-\infty$, $a_2=0$, $\epsilon_2=-$, 
$\sigma_1=+$, which gives $\int^{G(0-)}_{G(-\infty)+}f=\int_{1+}^{\infty}f
=F(\infty)-F(1+)=\int_{(-\infty,0)}(f\circ G)G'$.

\begin{theorem}[Translations] {\rm (a)} $\alexr$ is invariant under translation, i.e.,
$f\in\alexr$ if and only if $\tau_t f\in\alexr$ for  all $t\in\R$.
{\rm (b)} $\|\tau_t f\|=\|f\|$ for all $f\in\alexr$ and all $t\in\R$.
\end{theorem}

\bigskip
\noindent
{\bf Proof:} 
(a) Let $f\in\alexr$.  Then $f=F'$ for $F\in\balexr$.  For $\phi\in\D$ we have
\begin{eqnarray*}
\langle(\tau_tF)',\phi\rangle & = & -\langle\tau_tF,\phi'\rangle
=-\langle F,\tau_{-t}\phi'\rangle =-\langle F,(\tau_{-t}\phi)'\rangle\\
 & = & \langle F',\tau_{-t}\phi\rangle = \langle \tau_tF',\phi\rangle
  =  \langle \tau_tf,\phi\rangle.
\end{eqnarray*}
If $f\in\Dp$ such that $\tau_tf\in\alexr$, reverse  the above steps.

(b) Note that $\|\tau_tf\| = \sup_{x\in\R}|\tau_tF(x)|=\sup_{x\in\R}|F(x-t)|
=\|F\|_\infty=\|f\|$.  \qed\\

We have  {\it continuity in norm} if for $f\in\alexr$ we have
$\|f-\tau_xf\|\to 0$ as $x\to 0$.  But this is not true in $\alexr$.
For example,
$\|H_1'-\tau_x H_1'\|=1$ if $x\not=0$.

\section{Integration by parts}\label{Sectionintegrationbyparts}
An integration by parts
formula is obtained  using the Henstock--Stieltjes
integral.  (See the Introduction.)
This allows us to prove versions of the H\"older
inequality and Taylor's theorem.

The integration by parts  formula in $\alexr$ follows  from integration
by parts for the Henstock--Stieltjes integral \cite[p.~199]{mcleod}.
There it is proved that
$\int_{-\infty}^\infty F\,dg$ and $\intinf g\,dF$ exist when  one of $F$ and $g$ is
regulated and one is of bounded variation.
See also
\cite{tvrdy1989} and \cite{tvrdy2002} 
where various properties of these integrals are
established.

We first  need to define the product
of $f\in\alexr$ and  $g\in\bv$.

\begin{proposition}\label{Psi}
For $f\in\alexr$ and $g\in\bv$,  let $\{c_n\}$ contain all
$t\in\R$ such that
both $F$  and  $g$  are not right continuous
at  $t$. Define
$\Psi(x)=F(x)g(x)-\int_{-\infty}^xF\,dg-\sum_{c_n<x}\left[
F(c_n)-F(c_n+)\right]\left[g(c_n)-g(c_n+)\right]$.
Then $\Psi\in\balexr$.
The sum is over all $n\in\N$ such that $c_n<x$.
The integral and series defining  $\Psi$  converge
absolutely.

\end{proposition}

\bigskip
\noindent
{\bf Proof:} There is $M\in\R$ such  that $|g|\leq M$ and $Vg\leq M$.
Let $x\in\R$. Then 
\begin{eqnarray*}
|\Psi(x)| & \leq & |F(x)|M + \|F\chi_{(-\infty,x]}\|_\infty Vg
+2\|F\chi_{(-\infty,x+1]}\|_\infty Vg\\
 & \to & 0 \quad\text{as } x\to-\infty.
\end{eqnarray*}
This also  shows that the integral and series defining  $\Psi$  converge
absolutely.
Let $y<x$.  Then
\begin{eqnarray*}
\Psi(y)-\Psi(x) & = & \left[F(y)-F(x)\right]g(y)+\int_{y}^x \left[F(t)-F(x)\right]
dg(t)\\
 & & \quad -
\lsum_{c_n\in[y,x)}\left[
F(c_n)-F(c_n+)\right]\left[g(c_n)-g(c_n+)\right]
\end{eqnarray*}
so that, using the uniform left continuity of $F$ 
(Theorem~\ref{theorembalexr}),
\begin{eqnarray*}
\lim_{y\to x^-}|\Psi(y)-\Psi(x)| & \leq & \lim_{y\to x^-}|F(y)-F(x)|M+ 
2\lim_{y\to x^-}\sup_{s,t\in[y,x]}|F(s)-F(t)| Vg\\
 & = & 0.
\end{eqnarray*}
Therefore, $\Psi$ is left continuous.  Similarly, using the
uniform right regularity of $F$ we see that $\Psi$ has a right limit
at each point.  Letting $x,y\to\infty$ in
the above inequality shows $\Psi(\infty)$ exists.\qed\\

If $F$ is taken as regulated but not left continuous then there is an additional
term in $\Psi$ involving $F(c_n)-F(c_n-)$.  See \cite[p.~199]{mcleod}.

For an arbitrary distribution $T\in\Dp$ we  have the product
$T\psi$ defined for all $\psi\in C^\infty(\R)$ by
$\langle T\psi,\phi\rangle=\langle T,\psi\phi\rangle$ for
$\phi\in\D$.  Distributions in $\alexr$ can be multiplied
by functions of bounded variation.
\begin{defn}[Product]\label{productdefn}  With $\Psi$ as in 
Proposition~\ref{Psi} and $\phi\in\D$, the product of
$f\in\alexr$ and $g\in\bv$ is defined by $\langle fg, \phi\rangle
=\langle \Psi',\phi\rangle=-\langle \Psi,\phi'\rangle$.
\end{defn}
This defines
$fg\in\alexr$ since $\Psi\in\balexr$.  Each of
the three terms in $\Psi$ is regulated so the product $\Psi(t)\phi(t)$
is Riemann integrable.

\begin{defn}[Integration by parts]\label{integrationbyparts}
Let $f\in\alexr$ and $g\in\bv$ and use  the notation of
Proposition~\ref{Psi}.  Define the integral of $fg$ by
\begin{eqnarray*}
\int_{-\infty}^\infty fg & = & \int_{-\infty}^\infty g\,dF\\
 & = & F(\infty)g(\infty)-\int_{-\infty}^\infty 
F\,dg - \sum_{n\in\N}\left[
F(c_n)-F(c_n+)\right]\left[g(c_n)-g(c_n+)\right].
\end{eqnarray*}
\end{defn}

Notice that if $F$ is continuous or if $g$ is right continuous
then the sum in the integration by parts formula vanishes and
we recover the more familiar formula $\intinf fg= F(\infty)g(\infty)
-\intinf F\,dg$.  Note also that we have
defined the integration by parts formula to agree with the
Stieltjes  integral but  we have  no way of proving the
formula.  However, when $F$ is appropriately smooth
it reduces to the usual formula for Lebesgue ($F\in AC$), Henstock--Kurzweil
($F\in ACG*$)
and wide Denjoy integrals ($F\in ACG$).
Density arguments show we have the correct formula in $\alexr$.
Since step functions are dense in the regulated functions 
\cite[\S7.13]{mcleod}, given
$f\in\alexr$ there is a sequence of step functions $\{F_n\}\in\balexr$
such that $\|F_n-F\|_\infty\to 0$.  Definition~\ref{integrationbyparts}
certainly holds for $f_n=F_n'$ and $g\in\bv$.  To see this it suffices
to prove the formula for $f=\delta$, $F=H_1$ and $g\in\bv$.  We have
$F(\infty)g(\infty)=g(\infty)$.  To evaluate $\intinf F\,dg$, take a gauge
$\gamma$ that forces $0$ to be a tag.  If $\Partition=\{([x_{n-1},x_n],z_n)\}_{n=1}^N$ is $\gamma$-fine then
\begin{eqnarray*}
\sum_{n=1}^N H_1(z_n)[g(x_n)-g(x_{n-1})] & = & 
\sum_{z_n>0}[g(x_n)-g(x_{n-1})]\\
 & = & g(\infty)-g(z)
\end{eqnarray*}
where $z$ is the smallest positive tag in $\Partition$.  We can take
$\gamma$ so that $z$ is as close to $0$ as we like.
Therefore, $\int_{-\infty}^\infty F\,dg=g(\infty)-g(0+)$.  And,
$-\sum[H_1(c_n)-H_1(c_n+)][g(c_n)-g(c_n+)]=g(0)-g(0+)$.  Hence,
$F(\infty)g(\infty)-\int_{-\infty}^\infty 
F\,dg - \sum_{n\in\N}\left[
F(c_n)-F(c_n+)\right]\left[g(c_n)-g(c_n+)\right] =g(0)$.  And,
we also have $\sum_{n=1}^N g(z_n)[H_1(x_n)-H_1(x_{n-1})]=g(0)$
so that $\intinf g\,dF=g(0)$.
   The H\"older inequality
(Theorem~\ref{holder} below) then gives $|\intinf(f_n-f)g|\leq
\|F_n-F\|_\infty\|g\|_{\bv}\to 0$ as $n\to\infty$.  This justifies
the integration by parts formula.

The calculation above shows that $\intinf \delta g =g(0)$ for each function
$g$ that has a right limit at $0$ and a limit at $\infty$.  
Thus, the integration by parts formula
is in accordance with the action of $\delta$ as a measure.
For example, let $a\in\R$ and define
\begin{equation}
g_a(x)=\left\{\begin{array}{cl}
0, & x<0\\
a, & x=0\\
1, & x>0.
\end{array}
\right.\label{onepoint}
\end{equation}
Then for $\phi\in\D$ we have $\langle\delta g_a,\phi\rangle
=\langle \delta, g_a\phi\rangle =a\phi(0)$.  Putting $a=0$ gives
$\langle \delta H_1,\phi\rangle=0$ so $\delta H_1\in\alexr$ and
$\delta H_1=0$; $a=1$ gives $\langle \delta H_2,\phi\rangle=\phi(0)$ so 
$\delta H_2\in\alexr$ and $\delta H_2=\delta$.  See \cite{raju}
for references to other methods of multiplying the Dirac and 
Heaviside distributions.
We also see that
changing $g$ at even one point can affect the value of
the integral of $fg$, i.e., for $f=\delta$ the integral depends on the value of $g(0)$.  And, defining a function $F$ to be the right side of \eqref{onepoint}
we see that the integration by parts formula does not depend on our convention
of using left continuous primitives.  For such $F$ and any $g\in\bv$, both
right sides of Definition~\ref{integrationbyparts} give zero, provided
we use the more general formula \cite[p.~199]{mcleod} 
that allows discontinuities from the
left and right.

The integration  by parts formula leads to a version of the
H\"older inequality.  Note that $\bv$ is a Banach space
under  the  norm $\|g\|_{\bv}=\|g\|_\infty+Vg$.

\begin{theorem}[H\"older]\label{holder}
Let $f\in\alexr$ and $g\in\bv$.  Then $\left|\int_{-\infty}^\infty
fg\right|\leq\\
\left|\int_{-\infty}^\infty
f\right||g(\infty)|+\|f\|Vg\leq \|f\| \|g\|_{\bv}$. 
The inequality is sharp in the sense that if 
$\left|\int_{-\infty}^\infty
fg\right|\leq \|f\|\left(\alpha|g(\infty)|+\beta Vg\right)$ for
all $f\in\alexr$ and all $g\in\balexr$ then $\alpha, \beta\geq 1$.
For each 
$-\infty\leq a\leq \infty$ there is the inequality 
$\left|\int_{-\infty}^\infty
fg\right|\leq \|f\|\left(|g(a)|+2Vg\right)$.
\end{theorem}

\bigskip
\noindent
{\bf Proof:}
Use the fact that $\int_{-\infty}^\infty fg=\int_{-\infty}^\infty g\,dF$.
Given $\epsilon>0$ there is a partition
$\{(z_n,[x_{n-1},x_n])\}_{n=1}^N$ so that
$|\intinf g\,dF-\sum_{n=1}^N g(z_n)[F(x_n)-F(x_{n-1})]|<\epsilon$.  Since
$F(x_0)=F(-\infty)=0$,
\begin{eqnarray*}
\left|\int_{-\infty}^\infty fg\right|  & \leq & \epsilon
+ \left|\sum_{n=1}^N g(z_n)[F(x_n)-F(x_{n-1})]\right|\\
 & = & \epsilon +  \left|\sum_{n=1}^N g(z_n)F(x_n)-
\sum_{n=1}^{N-1} g(z_{n+1})F(x_n)
\right|\\
 & = & \epsilon +  \left|F(\infty)g(\infty)-\sum_{n=1}^{N-1}
F(x_n)\left[g(z_{n+1})-g(z_n)\right]\right|\\
 & \leq & \epsilon + \left|\int_{-\infty}^\infty f\right| |g(\infty)|
+\|F\|_\infty Vg\\
 & \leq & \epsilon +\|f\| \|g\|_{\bv}.
\end{eqnarray*}
The final estimate follows upon  noting  that $|g(\infty)|\leq
|g(a)|+|g(\infty)-g(a)|\leq |g(a)|+Vg$.

We can see the estimate is sharp by letting $F(x)=(P/\pi)(\pi/2+
\arctan(x))$ and
$g(x)=Q$ for $x\leq a$ and $g(x)=R$ for $x>a$, where $P>0$ and
$Q>R>0$.
Then $\int_{-\infty}^\infty fg=PR+(Q-R)F(a)$.  As $a\to\infty$ we
see this approaches $F(\infty)g(\infty)+\|f\|Vg$.\qed\\

For a proof using the Henstock--Stieltjes integral see 
\cite[Theorem~2.8]{tvrdy1989} and \cite{tvrdy2002}.

Integration by parts and the fundamental theorem can be used
to prove a version of Taylor's theorem.
\begin{theorem}[Taylor]\label{taylor}
Let $f\fn[a,\infty)\to\R$. Let $n\in\N$.   If $f\in C^{n-1}([a,\infty))$
so that $f^{(n)}$ is regulated and
right continuous on $[a,\infty)$
then for all $x\in(a,\infty)$ we have $f(x) =P_n(x) + R_n(x)$ where
$$
P_n(x)=\lsum_{k=0}^n\frac{f^{(k)}(a)(x-a)^k}{k!}
\quad \text{ and } \quad
R_n(x)=\frac{1}{n!}\int_{(a,x]} f^{(n+1)}(t)(x-t)^{n}\,dt.
$$
We have the estimates
$
|R_n(x)|\leq \sup_{a\leq t< x}
|f^{(n)}(t)-f^{(n)}(a)| (x-a)^n/n!$ for $x\in(a,\infty)$ and
$\|R_n\chi_{(a,b)}\|\leq \|R_n\chi_{(a,b)}\|_1\leq (b-a)^{n+1}\sup_{a\leq t< b}
|f^{(n)}(t)-f^{(n)}(a)|/(n+1)!$.
\end{theorem}
Integration by parts gives an induction proof.
The remainder exists because
the function $t\mapsto (x-t)^n$ is in $\bv$ for each $x$.
The estimates on the remainder follow from the H\"older inequality
(Theorem~\ref{holder}).
Note that $R_n(x)=o((x-a)^n)$ as $x\to a+$.
If $f^{(n)}$ is left continuous on $(-\infty,a]$ then we expand $f$ in
powers of $a-x$ and $R_n(x)=o((a-x)^n)$ as $x\to a-$.
Usual versions
of Taylor's theorem require $f^{(n+1)}$ to be integrable.  For the Lebesgue
integral this means taking $f^{(n)}$ to be absolutely continuous.  
Here we only require $f^{(n)}$ to be regulated.
The case $n=0$ corresponds to Theorem~\ref{fundamental}.

\section{Norms and dual space}
Multipliers are 
those functions
$g$ for which $fg$ is integrable for all integrable $f$.  
In this section we consider some equivalent norms on $\alexr$ and then
show that the space
of multipliers of $\alexr$
and the dual space of $\alexr$ are both given by $\bv$.  
\begin{theorem}[Equivalent norms]\label{eqnorms}
{\rm (a)} The following norms are equivalent to $\|\cdot\|$ in $\alexr$.
For
$f\in\alexr$, define $\|f\|'=
\sup_I|\int_If|$, where the supremum is taken over all finite intervals 
$I\subset\R$;
$\|f\|''=\sup_{g}\int fg$, where the supremum is taken over all
$g\in\bv$ such that $\|g\|_\infty\leq 1$ and $Vg\leq 1$.
{\rm (b)} Let $g\in\bv$ and be normalised so that $g(x)=(1-\lambda) g(x-)+\lambda g(x+)$ 
for fixed $0\leq\lambda\leq 1$ and all $x\in\R$.   Then for $f\in\alexr$ we have
$\left|\int_{-\infty}^\infty
fg\right|\leq
\left|\int_{-\infty}^\infty
f\right|\inf|g|+\|f\|'\,Vg$.
{\rm (c)} Let $a\in\Rbar$.  The norms $\|g\|'_{\bv}=|g(a)|+Vg$ and
$\|g\|_\bv=\|g\|_\infty+Vg$ are equivalent on $\bv$.
\end{theorem}

\bigskip
\noindent
{\bf Proof:} (a) Note that  $\|f\|\leq \|f\|'$.
And, we have $\|f\|'=\sup_{a<b}\left|\int_{(a,b)} f\right|=
\sup_{a<b}\left| F(b-)-F(a+)\right|\leq 2\|f\|$. Similarly for other
types of intervals.
Hence, $\|\cdot\|$ and $\|\cdot\|'$ are equivalent.  Let $g\in\bv$
such that $\|g\|_\infty\leq 1$ and $Vg\leq 1$. 
By the H\"older
inequality (Theorem~\ref{holder}),
$$
\left|\intinf fg\right| \leq \|f\|\left[|g(\infty)|+Vg\right]\leq 2\|f\|.
$$
And,
$$
\|f\|''\geq \max\left(\sup_{x\in\R}\intinf f\chi_{(-\infty,x]},
-\sup_{x\in\R}\intinf f\chi_{(-\infty,x]}\right).
$$
It follows that $\frac{1}{2}\|f\|''\leq \|f\|\leq \|f\|''$.  

(b) The alternative H\"older inequality is proved as Lemma~24
in \cite{talvilafourier}. 

(c)
Clearly, $\|g\|'_\bv\leq \|g\|_\bv$. Let $x\in\R$.  The inequality
$|g(x)|\leq |g(a)|+
 |g(a)-g(x)|\leq |g(a)|+Vg$ shows $\|g\|_\bv\leq 2\|g\|'_\bv$.  \qed

The H\"older inequality can be reformulated in any of these
equivalent norms.

For the
Henstock--Kurzweil and continuous primitive integral \cite{talviladenjoy}
the multipliers
and dual space
are the functions of essential  bounded  variation.
(See Example~\ref{examples}(d) for the definition.)  
For  $\alexr$ the  multipliers and dual space are the  functions
of  bounded variation.
\begin{theorem}\label{dualalexr}
The set of multipliers for $\alexr$ is $\bv$.
\end{theorem}

\bigskip
\noindent
{\bf Proof:}
The multipliers are defined in Definition~\ref{integrationbyparts}.
Hence,  every function of bounded  variation is a multiplier.
In order for a function $g$ to be a multiplier the integral
$\intinf g\,dF$ must exist for every $F\in\balexr$.  Taking
$F$ to be a step function $F=\sum\sigma_n\chi_{(a_n,b_n]}$ for
disjoint intervals $\{(a_n,b_n)\}$  and $\sigma_n\in\R$, we
see that $\intinf g\,dF=\sum\sigma_n[g(a_n)-g(b_n)]$.  Taking
$\sigma_n={\rm sgn}(g(a_n)-g(b_n))$ shows $g\in\bv$.
\qed

If $\{f_n\}$ is a sequence in $\alexr$ such that $\|f_n\|\to 0$ then  the
H\"older inequality shows that $\int_{-\infty}^\infty f_ng\to 0$ for
each $g\in\bv$.  Hence, for each fixed $g\in\bv$, $T_g(f):=\int_{-\infty}
^\infty fg$ defines a  continuous linear functional on $\alexr$.
Hence, $\alexr^\ast\supset\bv$.  In fact, all continuous linear functionals
on $\alexr$ are of this form, i.e., $\alexr^\ast=\bv$.  We can prove
this by using the representation of the dual of the space of regulated
functions.

Various authors have used
different specialised integrals to represent the continuous linear
functionals on regulated functions.  See  Kaltenborn \cite{kaltenborn} 
(Dushnik interior integral) for compact intervals (also \cite{honig}), 
Hildebrandt \cite{hildebrandt1934} (refinement
or Young integral) for $\R$, Tvrd\'y \cite{tvrdy1989}, \cite{tvrdy1996},
\cite{tvrdy2002}
(Henstock--Stieltjes integral, where it is called the Perron--Stieltjes 
integral).
Tvrd\'y gives a representation for such a functional acting on
regulated function $F$ on compact interval $[a,b]$ as
$T(F) = qF(a)+\int_a^b g\,dF$ for some function $g\in\bv$ and $q\in\R$.
See also \cite{schwabik1992}.
Extension to regulated functions on $\R$ follows by replacing
$a$ with $-\infty$ and $b$ with $\infty$, using our definition
of the Henstock--Stieltjes integral 
(Section~\ref{Sectionintegrationbyparts})
and  compactification of $\R$ (Remark~\ref{compactification}).
For $F\in\balexr$, the functional
then becomes $T_g(F)=\intinf g\,dF=\intinf F'g$.  The connection between the
Dushnik interior and Young integrals is given in \cite{hildebrandt1966}.
Equality of Young and Henstock--Stieltjes integrals for one function
regulated and one of bounded variation is established in \cite{schwabik1973}.

\begin{theorem}\label{dual}
The dual space of $\alexr$ is $\bv$ ($\alexr^\ast =\bv$).
\end{theorem}
\bigskip
\noindent
{\bf Proof:}  Let $\psi\fn\alexr\to\balexr$ be  given by $\psi(f)=F$.
Then $\psi^{-1}\fn\balexr\to\alexr$ is given by $\psi^{-1}(F)=F'$.
Let $\{f_n\}\subset\alexr$ such that $\|f_n\|\to0$.
Then $\|F_n\|_\infty\to0$.  If $T\in\alexr^\ast$ then $T(f_n)=
T(\psi^{-1}(F_n))\to0$.  Hence, $T\circ\psi^{-1}\in\balexr^\ast$.  
Using  the result of the previous paragraph, we have
$\balexr^\ast =\bv$.    There exists
$g\in\bv$ such that $T\circ\psi^{-1}(F_n)=\intinf F_n'\,dg=\intinf f_ng$.
Hence, $T(f_n)=\intinf f_ng$.
\qed

The integration by parts formula also shows that $\langle f,g\rangle
=\int_{-\infty}^\infty fg$ for all $f\in\alexr$ and all $g\in\bv$
so that we could use  integration by parts as a starting point
to define the integral
as a continuous linear functional on $\bv$.

In the space of Henstock--Kurzweil integrable functions we identify
functions almost everywhere so the dual of this space is
$\ebv$ rather than $\bv$, i.e., if $T$ is a continuous  linear functional
on the space of
Henstock--Kurzweil integrable functions then  there exists a function
$g\in\bv$ such  that $T(f)=\intinf fg$ for each Henstock--Kurzweil
integrable function  $f$.  The integral is that of  Henstock--Kurzweil.
Changing $g$ on a set of measure zero does not affect the value of
this integral  so the dual space is $\ebv$.

In $\alexr$ we do not have this equivalence relation so 
the dual of $\alexr$ is $\bv$ and not $\ebv$.  Similarly, for
no normalisation in $\bv$ (see Remark~\ref{equivalencerelation}) is the
dual of $\alexr$ equal to  functions of normalised bounded variation.
To see this, note that the function $g=\chi_{\{0\}}$ is not
equivalent to $0$ since $\int_{-\infty}^\infty fg=F(0+)-F(0-)$.
But every normalisation makes $g=0$.

No concrete description of $\bv^\ast$ seems to
be known. 
But note that if $\{g_n\}\subset \bv$ such that $\|g_n\|_\bv\to 0$ then
$\int_{-\infty}^\infty f g_n\to 0$ for each $f\in\bv$.  Hence,
$T_f(g)=\int_{-\infty}^\infty fg$ defines a continuous linear functional
on $\bv$.  The H\"older inequality shows that for 
each regulated function $F$ the linear functional
\begin{eqnarray}
T_F(g) & = & \intinf F\,dg =- \int_{-\infty}^\infty F' g
+ F(\infty)g(\infty) -F(-\infty)g(-\infty)\notag\\
 & & \quad- \sum_{n\in\N}\left[
F(c_n)-F(c_n+)\right]\left[g(c_n)-g(c_n+)\right]\label{linearfunctional}\\
 & & \quad +\sum_{n\in\N}\left[
F(c_n)-F(c_n-)\right]\left[g(c_n)-g(c_n-)\right]\notag
\end{eqnarray}
is  in $\bv^\ast$, i.e., if $\|g_n\|_{\bv}\to 0$ then
$T_F(g_n)\to 0$ in $\R$.  Hence, $\bv^\ast$ contains the space of regulated
functions.  If we let $F=\chi_{\{0\}}$
then  $T_F(g)=\intinf F\,dg=g(0+)-g(0-)$ and $|T_F(g)|\leq Vg$ so
$T_F\in\bv^\ast$ but as an element of $\alexr$, $F'=0$.  Hence,
$\bv^\ast\supsetneq \alexr$.
And, consider
the following example.  Let $S=\{1/n\mid n\in\N\}$, $F=\chi_S$ and
define $U_F\fn\bv\to\R$ by $U_F(g)=\intinf F\,dg$.  Then $F$ is not
of bounded variation and since $\lim_{x\to 0+}F(x)$ does not exist, 
$F$ is not regulated.  But, for $g\in\bv$, $U_F(g)=\sum_{n=1}^\infty
[g(n^{-1}+)-g(n^{-1}-)]$.  This can be seen by taking a gauge that
forces $0$ to be a tag and forces
$n^{-1}$ to be  a tag for some $N_0$ and all $1\leq n\leq N_0$.  We  then
have $|U_F(g)|\leq Vg$.  This shows that
$U_F\in\bv^\ast$.   Hence, $\bv^\ast$ properly
contains the space of regulated functions.  More precisely, the
space of regulated functions is identified with finitely additive
measures defined by $\mu([a,b])=F(b+)-F(a-)$  for regulated function
$F$.  Similarly for other intervals.  These measures are defined on the algebra generated by intervals.

Hildebrandt
\cite{hildebrandt1966} and Aye and Lee \cite{aye} have given
explicit representation of the dual of $\bv^\ast$ in the topology of uniform
bounded variation with uniform convergence.
A sequence $\{g_n\}\subset\bv$ converges to $0$ in this sense if
$\|g_n\|_\infty\to 0$ and there is $M\in\R$ so that for all $n\in\N$,
$Vg_n\leq M$.  These authors show that the dual
of $\bv$ in this topology contains only (pairs of) regulated functions.  This
dual must then be a proper subset of $\bv^\ast$ since
$U_F$ from the preceding paragraph 
is not continuous in the topology of uniform
bounded variation with uniform convergence.
For example, define the piecewise linear functions
$g_n(x)=n^{-1}\sum_{m=1}^n(1-m+m^2x)\chi_{[m^{-1},m^{-1}+m^{-2}]}(x)$.
Then $\|g_n\|_\infty =1/n$ and $Vg_n=1$.  Hence, $g_n\to 0$ in the 
topology of \cite{hildebrandt1966} and \cite{aye}.  But, 
$$
U_F(g_n)  =  n^{-1}\sum_{m=1}^n F(m^{-1})[g(m^{-1}+)-g(m^{-1}-)]
 =   1 \not\to 0.
$$

Mauldin (\cite{mauldin} and references therein) and Hildebrandt 
\cite{hildebrandt1938} have 
given representations of $\bv^\ast$ in terms of abstract
integrals.

\section{$\bv$-module}
In Definition~\ref{productdefn}
we have a
product defined from $\alexr\times\bv$ onto $\alexr$.  
It has distributive, commutative
and associative properties that make  $\alexr$  into a Banach $\bv$-module.
See \cite{dales} for the   definition.
Properties of the integral 
of $fg$ then follow from properties of the product.

\begin{theorem}[Products]\label{producttheorem}
Let $f,f_1,f_2\in\alexr$; $g,g_1,g_2\in\bv$; $k\in\R$.
The
product has  the  following properties.
{\rm (a)} Distributive. $(f_1+f_2)g=f_1g+f_2g$, $f(g_1+g_2)=fg_1+fg_2$.
{\rm (b)} Homogeneous.  $(kf)g=f(kg)=k(fg)$.
{\rm (c)} Commutative. $f(g_1g_2)=f(g_2g_1)$.
{\rm (d)} Compatible with distribution product. 
$\langle fg,\phi\rangle= \langle f, g\phi\rangle$ for all $\phi\in\D$.
{\rm (e)} Associative. $(fg_1)g_2=f(g_1g_2)$.
{\rm (f)} Zero divisors.  There are $f\not=0$ and $g\not=0$ such that $fg=0$.
{\rm (g)} Compatible with pointwise product. If $f$ and $g$ are functions
that are continuous at $a\in\R$ then $\langle fg,\phi_n\rangle \to f(a)g(a)$
for any $\delta$-sequence supported at $\{a\}$.
\end{theorem}

\bigskip
\noindent
{\bf Proof:} Properties (a), (b) and (c) 
follow immediately from the  definition.

Notice we can include terms in the sum (also labeled $c_n$ but
not necessarily points were $F$ and $g$ are simultaneously 
discontinuous from the  right) so that $\sup c_n=\infty$.
To prove (d), let $f\in\alexr$, $g\in\bv$  and $\phi\in\D$.  Then
\begin{eqnarray*}
\langle fg,\phi\rangle & = & -\int_{-\infty}^\infty F(x)g(x)\phi'(x)\,dx
+ \int_{-\infty}^\infty \int_{-\infty}^x F(t)\,dg(t)\,\phi'(x)\,dx\\
 & & \quad+\int_{-\infty}^\infty \lsum_{c_n<x}\left[
F(c_n)-F(c_n+)\right]\left[g(c_n)-g(c_n+)\right]\phi'(x)\,dx.
\end{eqnarray*}
Note that $|\int_{-\infty}^\infty \int_{-\infty}^x F(t)\,dg(t)\,\phi'(x)\,dx|\leq \|F\|_\infty
Vg \|\phi'\|_1$ so by the Tonelli  and Fubini theorems,
$\int_{-\infty}^\infty \int_{-\infty}^x F(t)\,dg(t)\,\phi'(x)\,dx
=-\int_{-\infty}^\infty
F(t)\phi(t)\,dg(t)$.  And,
$|\int_{-\infty}^\infty \sum_{c_n<x}\left[
F(c_n)-F(c_n+)\right]\left[g(c_n)-g(c_n+)\right]\phi'(x)\,dx|\leq 2\|F\|_\infty
Vg\, \|\phi'\|_1$ so  again,
\begin{eqnarray*}
\lefteqn{
\int_{-\infty}^\infty \sum_{c_n<x}\left[
F(c_n)-F(c_n+)\right]\left[g(c_n)-g(c_n+)\right]\phi'(x)\,dx}\\
 & = &
\sum_{n=1}^\infty \left[
F(c_n)-F(c_n+)\right]\left[g(c_n)-g(c_n+)\right]
\int_{x>c_n}\phi'(x)\,dx\\
 & = & -\sum_{n=1}^\infty \left[
F(c_n)-F(c_n+)\right]\left[g(c_n)-g(c_n+)\right]\phi(c_n).
\end{eqnarray*}
Therefore, 
\begin{eqnarray}
\langle fg,\phi\rangle & = & -\int_{-\infty}^\infty F(t)g(t)\phi'(t)\,
dt-\int_{-\infty}^\infty F(t)\phi(t)\,dg(t)\notag\\
 & & \quad -
\sum_{n=1}^\infty \left[
F(c_n)-F(c_n+)\right]\left[g(c_n)-g(c_n+)\right]\phi(c_n).\label{fgphi}
\end{eqnarray}
And, $g\phi\in\bv$ with compact support, so using the continuity of $\phi$,
\begin{eqnarray*}
\lefteqn{
\langle f,g\phi\rangle=\int_{-\infty}^\infty f(g\phi)}\notag\\ & = & 
-\int_{-\infty}^\infty F\,d(g\phi)
-\sum \left[
F(c_n)-F(c_n+)\right]\left[g(c_n)\phi(c_n)-g(c_n+)\phi(c_n+)\right]\\
 & = & -\int_{-\infty}^\infty F\,d(g\phi)
-\sum \left[
F(c_n)-F(c_n+)\right]\left[g(c_n)-g(c_n+)\right]\phi(c_n).
\end{eqnarray*}
Now show that $\intinf F\,d(g\phi)=\intinf Fg\,d\phi+\intinf F\phi\,dg$.
Each of these integrals
exists because, in each case, one of the integrands and integrators
is of bounded variation and one is regulated.  
For $\epsilon>0$ there is then a tagged partition 
$\{(z_n,[x_{n-1},x_n])\}_{n=1}^N$
such that
$|S_N-\int_{-\infty}^\infty F\,d(g\phi)+
\int_{-\infty}^\infty Fg\,d\phi + \int_{-\infty}^\infty F\phi\,dg|<\epsilon$,
where 
\begin{eqnarray*}
S_N & = & \sum_{n=1}^N 
F(z_n)\left\{\left[g(x_n)\phi(x_n)-g(x_{n-1})\phi(x_{n-1})\right]\right. \\
 & & \left.\quad -g(z_n)\left[\phi(x_n)-\phi(x_{n-1})\right]
-\phi(z_n)\left[g(x_n)-g(x_{n-1})\right]\right\}.
\end{eqnarray*}
But, 
\begin{eqnarray*}
|S_N| & \leq & \sum_{n=1}^N 
|F(z_n)|\left\{|[g(x_n)-g(x_{n-1})||\phi(x_n)-\phi(z_{n})|\right.\\
 & & \quad\left.
+|g(x_{n-1})-g(z_n)||\phi(x_n)-\phi(x_{n-1})|\right\}.
\end{eqnarray*}
Since $\phi$ is uniformly continuous we can arrange the partition so
that the maximum of $|\phi(x_n)-\phi(t_n)|$ for
$t_n\in[x_{n-1},x_n]$ is less than $\epsilon$
for each $1\leq n\leq N$.  Then $|S_N|\leq 2\epsilon\|F\|_\infty Vg$.
Hence, $\int_{-\infty}^\infty F\,d(g\phi)=
\int_{-\infty}^\infty Fg\,d\phi + \int_{-\infty}^\infty F\phi\,dg$.
Now using \eqref{fgphi} we see that $\langle fg, \phi\rangle=
\langle f,g\phi\rangle$.

Associativity (e) then follows by writing
$\langle f(g_1g_2),\phi\rangle  =  \langle f,(g_1g_2)\phi\rangle =
\langle f,g_1(g_2\phi)\rangle
  =  \langle fg_1,g_2\phi\rangle = \langle (fg_1)g_2,\phi\rangle$.

To  prove (f), let $F\in\balexr$ and $g\in\bv$ be continuous
with disjoint  support.  Then  $F'g=0$.

{\it A $\delta$-sequence supported at $a$} is a sequence $\{\phi_n\}\subset
\D$ such that $\phi_n\geq 0$, $\intinf \phi_n=1$, $\supp(\phi_n)$
is an interval containing $a$ in its interior such that $\supp(\phi_n)
\to \{a\}$. For such a sequence, suppose $\supp(\phi_n)\subset[a-\delta,
a+\delta]$.  Then for (g),
\begin{eqnarray*}
\left|f(a)g(a)-\intinf fg\phi_n\right| & = & 
\left|\intinf\left[f(a)g(a)\phi_n- fg\phi_n
\right]\right|\\
 & \leq & \int_{a-\delta}^{a+\delta}\left|f(a)g(a)-f(x)g(x)\right|\phi_n(x)
\,dx\\
 & \to & 0\quad\text{ using the continuity of $f$ and $g$}.\qed
\end{eqnarray*}

If $g\in C^\infty$ then we see the product reduces to the usual
product of a distribution and a smooth function (cf. the paragraph
preceding Definition~\ref{productdefn}).

Each result in Theorem~\ref{producttheorem} concerning a product
can be integrated.  For example, $\int_{-\infty}^\infty
f(gh)=\int_{-\infty}^\infty (fg)h$.  Taking $g$ to be 
the characteristic function of an interval
and integrating by parts recovers each of the four 
integrals defined in \eqref{intdefn1}-\eqref{intdefn4}: $\int_{-\infty}^\infty\
f\chi_{I}=\int_I f$ for any interval $I$.
Each of the integrals $\int_{-\infty}^\infty F\,d\chi_{I}$ and
$\int_{-\infty}^\infty \chi_I\,dF$ exists as a Henstock--Stieltjes
integral because we can take a gauge that forces endpoints of $I$ to
be tags.  If $F$ is not continuous at the endpoints of $I$ then
these integrals will not exist as Riemann--Stieltjes integrals.  

The usual pointwise product makes $\bv$ into an algebra with
unit $g=1$.  Our product on $\alexr\times \bv$ makes $\alexr$
into a (left) Banach $\bv$-module.   
\begin{theorem}[Banach $\bv$-module]
$\bv$ is a Banach algebra.
$\alexr$ is a Banach  $\bv$-module.
\end{theorem}

\bigskip
\noindent
{\bf Proof:}
For $g_1,g_2\in\bv$,  the
inequalities 
\begin{eqnarray*}
\|g_1g_2\|_\bv & = & \|g_1g_2\|_\infty+V(g_1g_2)\\
 & \leq & \|g_1\|_\infty\|g_2\|_\infty +\|g_1\|_\infty Vg_2 +Vg_1\|g_2\|_\infty
\\
 & \leq & \|g_1\|_\bv\,\|g_2\|_\bv
\end{eqnarray*}
show that $\bv$ is closed under multiplication.  It then  follows easily that
$\bv$ is a Banach algebra.

The second statement follows from (a), (b), (c) and (e) of 
Theorem~\ref{producttheorem} and the
inequality $\|fg\|\leq \|f\| \|g\|_\bv$, valid for  all $f\in\alexr$
and $g\in\bv$. 
\qed\\

Notice that $\bv$ is not a division ring since $\chi_{[0,1]}\not=0$
has no multiplicative inverse.  There are zero divisors.  For 
example, $\chi_{[0,1]}\chi_{[2,3]} =0$.

Notice that if $g_1, g_2\in\bv$ then $(g_1g_2)'=g_1'g_2+g_1g_2'\in\alexr$.
The product on the left is pointwise in $\bv$ while the products on
the right are as per Definition~\ref{productdefn}.  Hence, the distributional
derivative is a
{\it derivation} on the algebra $\bv$ into the  Banach $\bv$-module $\alexr$.
See \cite{dales}.

\section{Absolute integrability}\label{absoluteintegrability}
The primitives of an $L^1$ function are absolutely continuous
and hence are functions of bounded variation.  Whereas, if 
function $f$ is Henstock--Kurzweil or wide Denjoy integrable 
but $|f|$ is not integrable in this sense then the primitives
of $f$ are not of bounded variation.  We use this observation
to define absolute integrability in $\alexr$.
We also show that $L^1$ and the space of signed Radon measures
are embedded continuously in $\alexc$ and $\alexr$, respectively.
\begin{defn}[Absolute integrability, \nbv]
Define the functions of normalised bounded variation as 
$\nbv=\balexr\cap\bv$.
A distribution $f\in\alexr$ is {\it absolutely integrable}
if it has a primitive $F\in\nbv$.  Denote the space of
absolutely integrable distributions by $\anbv$.
\end{defn}

Hence, $\anbv$ is isometrically isomorphic to
the space of signed Radon measures under
the Alexiewicz norm.  For $f\in\anbv$, let its primitive
in $\nbv$ be $F$. As in Example~\ref{examples}(d)   
there is a unique signed Radon measure $\mu$
such that $F'=\mu$, i.e., $\langle F,\phi'\rangle=-\intinf \phi\,d\mu$
for all $\phi\in\D$.
And, if $\mu$ is a signed Radon measure then a function defined
by $F(x)=\mu((-\infty,x))$ is in $\nbv$.  The Alexiewicz 
norm of $\mu$ identified
with $f\in\anbv$ is  $\|\mu\|=\sup_{x\in\R}|\mu((-\infty,x))|$.

This then gives an alternative definition of the regulated
primitive integral.  It is the completion of the space of
signed Radon measures in the Alexiewicz norm.
Integration in
$\anbv$ is thus Lebesgue integration.

Denote the space of signed Radon measures by $\M$.  A norm
is given by $\|\mu\|_\M=|\mu|(\R)=\mu^+(\R)+\mu^-(\R)$, which
is the total variation of $\mu$. 
\begin{theorem}\label{theoremnbv}
{\rm (a)} $\nbv$ is a Banach subspace of $\balexr$ under the
norm $\|g\|_{\bv}=\|g\|_\infty +Vg$. 
{\rm (b)} For each $a\in\Rbar$ the norms $\|g\|_{\bv}$
and $\|g\|_{\bv a} :=|g(a)|+Vg$ are equivalent.
{\rm (c)} $\nbv$ is not a Banach space under $\|\cdot\|_\infty$.
{\rm (d)} $\nbv$ is dense in $\balexr$.  The completion of $\nbv$
in $\|\cdot\|_\infty$ is  $\balexr$.
{\rm (e)} $\anbv$ is a Banach subspace of $\alexr$.
{\rm (f)} $\alexr$ is the completion of the space of signed
Radon measures in the Alexiewicz
norm.
{\rm (g)}
The embeddings $L^1\hookrightarrow\alexc$ and
$\M\hookrightarrow\alexr$ are continuous.
\end{theorem}

\bigskip
\noindent
{\bf Proof:} (a) It is a classical result that functions of 
bounded variation form a Banach space.  For example, see
\cite{kannan}.  The case of $\nbv$ is similar, as in the proof of
Theorem~\ref{theorembalexr}.  (b)  See Theorem~\ref{eqnorms}(c).
(c) Let $g(x)=x\sin(x^{-2})$ for $x>0$ and $g(x)=0$ for $x\leq 0$.
Then $g\in C^0\setminus\bv$.  Let $g_n=[1-\chi_{[0,(n\pi)^{-1/2}]}]g$.
Each $g_n\in\nbv$.  And, $\|g_n-g\|_\infty\leq (n\pi)^{-1/2}\to 0$.
In $\|\cdot\|_\infty$ the sequence $\{g_n\}$ converges to $g\notin\bv$.
(d) Let $F\in\balexr$ and let $\epsilon>0$ be given.  There exists
$M>0$ such that $|F(x)|<\epsilon$ for all $x\leq-M$
and $|F(x)-F(\infty)|<\epsilon$ for all $x\geq M$.  For each
$x\in[-M,M]$ there is $\delta_x>0$ such that if $y\in(x-\delta_x,x]$
then $|F(y)-F(x)|<\epsilon$ and if $y\in(x,x+\delta_x)$ then
$|F(y)-F(x+)|<\epsilon$.  Let $I_x=(x-\delta_x,x+\delta_x)$.  The
collection $\{I_x\}_{x\in[-M,M]}$ is an open cover of the compact
interval $[-M,M]$.  There is then a finite subcover, $\{I_x\}_{x\in J}$
for some finite set $J\in[-M,M]$.  We can then take open subintervals
$I_x'\subset I_x$ such that each point in $[-M,M]$ is in either one or
two of these intervals.  Then we can define $g(x)=0$ for $x\leq -M$,
$g(x)=F(\infty)$ for $x>M$ and $g$ is piecewise constant on each 
interval $I_x'$ such that $g\in\nbv$ and $\|g-F\|_\infty<\epsilon$.
Hence, $\nbv$ is dense in $\balexr$ and its completion is $\balexr$.
(e), (f) These follow from the isomorphism between $\alexr$ 
and $\balexr$ given by the integral. (g) For $\mu\in\M$ we
have 
$$
\|\mu\|  =  \lsup_{x\in\R}|\mu((-\infty,x))|
  \leq   \lsup_{x\in\R}\left[\mu^+((-\infty,x))+\mu^-((-\infty,x))\right]
  =  \|\mu\|_\M.
$$
If $f\in L^1$ then $\|f\|=\sup_{x\in\R}|\int_{-\infty}^xf|\leq \|f\|_1$.
\qed\\

In Proposition~\ref{finitelyaddivenbv} below it is shown that
distributions in $\alexr$ are finitely additive measures that
are finite when their primitives are of bounded variation.

Here is an
alternative way of defining functions of normalised bounded
variation.  Fix $0\leq \lambda\leq 1$.  For
$g\in\bv$ define $g_\lambda(x)=(1-\lambda) g(x-)+\lambda g(x+)$,
$g_\lambda(-\infty)=g(-\infty)$ and $g_\lambda(\infty)=g(\infty)$.
Define $\nbvl=\{g_\lambda\mid g\in\bv\}$.  Then $\nbvl$ is a
Banach space under $\|g\|_{\bv}=\|g\|_\infty+Vg$.  The connection
with the functions of essential bounded variation is the following.
As in Example~\ref{examples}(d) we have
$\ebv=\{g\in L^1_{loc}\mid {\rm ess\,var}\, g<\infty\}$.  This is a Banach
space under the norm 
$\|g\|_{{\cal EBV}}= {\rm ess\,sup}\, |g| + {\rm ess\,var}\, g$.
The space $\ebv$ consists of equivalence classes of functions
identified almost everywhere.  For each $g\in\ebv$ there is a unique
$g_\lambda\in\nbvl$ such that ${\rm ess\,sup}\, g=\|g_\lambda\|_\infty$
and ${\rm ess\,var}\, g=Vg_\lambda$.  For each $0\leq\lambda\leq 1$ the Banach
spaces $\nbvl$ and $\ebv$ are isometrically isomorphic.  These
spaces are distinct from $\bv$.  For example, $\chi_{\{0\}}$
is equivalent to $0$ in $\ebv$, its normalisation is $0$ in $\nbvl$
but  $V\chi_{\{0\}}=2$ in $\bv$.
Note that $\nbv$ is isometrically  isomorphic to the signed Radon
measures and to $\anbv$, whereas $\nbvl$ is isometrically isomorphic
to $\anbv\times \R$.  If $g\in\nbvl$ then its distributional derivative
is a Radon measure $\mu$ and $g(-\infty)\in\R$. 
For more on essential variation see \cite{ziemer}.

If $F\in\nbv$ then there are increasing functions of
normalised
bounded variation $G$ and $H$ such that $F=G-H$.  A
distribution $T$ is {\it positive} if $\langle T,\phi\rangle\geq 0$
for each $\phi\in\D$ with $\phi\geq 0$.  Suppose $\phi\geq 0$.
Let $[a,b]$ contain the support of $\phi$.  
By the second mean value theorem for integrals \cite[p.~211]{mcleod}
there is $\xi\in[a,b]$
such that
\begin{eqnarray*}
\langle G',\phi\rangle & = & -\int_a^b G\phi'
 =  -\left[G(a)\int_a^\xi \phi' + G(b)\int_\xi^b\phi'\right]\\
 & = & \left[G(b)-G(a)\right]\phi(\xi)
 \geq  0.
\end{eqnarray*}
Hence, $f\in\anbv$ if and only if it can be written as $f=G'-H'$
for $G,H\in\nbv$ with $G',H'\geq 0$.

In the next section we will introduce an ordering suitable for
all distributions in $\alexr$.

\section{Banach lattice}\label{lattice}
In $\balexr$ there is the partial order: $F\leq G$ if and only if
$F(x)\leq G(x)$ for all $x\in\R$.  Note that this order depends
on our choice that functions in $\balexr$ be left continuous.  Since
$\alexr$ is isomorphic to $\balexr$ it inherits this partial order.
For $f,g\in\alexr$ define $f\preceq g$ if and only if $F\leq G$, where
$F$ and $G$ are the respective primitives in $\balexr$.  This
order is not compatible with the usual order on distributions:
if $T,U\in\Dp$ then $T\geq U$ if and only if $\langle T-U,\phi\rangle
\geq 0$ for all $\phi\in\D$ such that $\phi\geq 0$.  Nor is it compatible
with pointwise ordering in the case of functions  in $\alexr$.  For
example, if $f(t)=H_1(t)\sin(t^2)$ then $F\geq 0$ so $f\succeq 0$ in
$\alexr$ but not pointwise.  And, $f$ is  not positive in the
distributional sense.  Note, however, that if $f\in\alexr$ is a measure
or  a nonnegative
function or distribution then $f\succeq 0$ in $\alexr$.

The importance of  this ordering is that it interacts with the Alexiewicz
norm so that $\alexr$ is a Banach lattice. If $\preceq$ is a binary
operation on set $S$ then it is a {\it partial order} if for all
$x,y,z\in S$ it is {\it reflexive} ($x\preceq x$), {\it antisymmetric}
($x\preceq y$ and $y\preceq x$ imply $x=y$) and {\it transitive} ($x\preceq y$
and $y\preceq z$ imply $x\preceq z$). 
If
$S$ is a Banach space with norm $\|\cdot\|_S$ and $\preceq$ is a partial
order on $S$  then $S$ is a {\it Banach lattice} if for all $x,y,z\in S$
\begin{enumerate}
\item
$x\vee y$ and $x\wedge y$ are in $S$.  The {\it join} is 
$x\vee y=
\sup\{x,y\}=w$ such that $x\preceq w$, $y\preceq w$ and if $x\preceq\tilde{w}$
and $y\preceq\tilde{w}$ then $w\preceq\tilde{w}$.
The {\it meet} is 
$x\wedge y=
\inf\{x,y\}=w$ such that $w\preceq x$, $w\preceq y$ and if $\tilde{w}\preceq x$
and $\tilde{w}\preceq y$ then $\tilde{w}\preceq w$.
\item
$x\preceq y$ implies $x+z\preceq y+z$.
\item
$x\preceq y$ implies $kx\preceq ky$ for all $k\in\R$ with $k\geq 0$.
\item
$|x|\preceq |y|$ implies $\|x\|_S\leq \|y\|_S$.
\end{enumerate}
If $x\preceq y$ we write $y\succeq x$.
We also define
$|x|=x\vee (-x)$, $x^+=x\vee 0$ and $x^-=(-x)\vee 0$. 
Then $x=x^+-x^-$ and $|x|=x^++x^-$.

We have absolute
integrability: if $f\in\alexr$ so  is $|f|$.
The lattice operations are defined for $F,G\in\balexr$ by
$(F\vee G)(x)=
\sup(F,G)(x)=\max(F(x),G(x))$.
And, $(F\wedge G)(x)=\inf(F,G)(x)=
\min(F(x),G(x))$.
\begin{theorem}[Banach lattice]
{\rm (a)} $\balexr$ is a Banach lattice.  
{\rm (b)}
For $f,g\in\alexr$, define $f\preceq g$ if $F\leq G$ in $\balexr$.
Then $\alexr$ is a Banach lattice isomorphic to $\balexr$.
{\rm (c)}  Let $F,G\in\balexr$.  Then $(F\vee G)'=F'\vee G'$,
$(F\wedge G)'=F'\wedge G'$, $|F'|=|F|'$,
$(F^+)'=(F')^+$, and $(F^-)'=(F')^-$.
{\rm (d)} If $f\in\alexr$ then $|f|\in\alexr$ with primitive $|F|\in\balexr$.  
For each interval $I\subset\R$ we have $|\int_If|\geq |\int_I|f||$.
For each $-\infty<x\leq \infty$ we have $|\int_{(-\infty,x)}f|
=\int_{(-\infty,x)}|f|$.  And, $\|\,|f|\,\|=\|f\|$, $\|f^{\pm}\|\leq \|f\|$.  
{\rm (e)} If $f\in\alexr$ then $f^{\pm}\in\alexr$ with respective
primitives $F^{\pm}\in\balexr$.  {\it Jordan decomposition}: $f=f^+ - f^-$.
And, $\int_If=\int_If^+-\int_If^-$ for every interval $I\subset\Rbar$.
{\rm (f)} $\alexr$ is {\it distributive}:
$f\wedge(g\vee h)=(f\wedge g)\vee(f\wedge h)$
and  $f\vee(g\wedge h)=(f\vee g)\wedge(f\vee h)$ for all $f,g,h\in\alexr$. 
{\rm (g)} $\alexr$ is {\it modular}: For all $f,g\in\alexr$, if
$f\preceq g$ then $f\vee(g\wedge h)=g\wedge(f\vee h)$ for all $h\in\alexr$.
{\rm (h)} Let $F$ and $G$ be regulated functions on $\R$ with real limits
at $\pm\infty$. Then 
\begin{eqnarray}
F'\preceq G' & \Longleftrightarrow & F(x-)-F(-\infty)\leq G(x-)-G(-\infty)
\quad\forall x\in\R\label{lattice1}\\
 &  \Longleftrightarrow & F(x+)-F(-\infty)\leq G(x+)-G(-\infty)
\quad\forall x\in\R.\label{lattice2}
\end{eqnarray}
\end{theorem}

\bigskip
\noindent
{\bf Proof:} (a) Let  $F,G\in\balexr$.  Define $\Phi=(F\vee G)$
and $\Psi=(F\wedge G)$.
We need to prove $\Phi,\Psi\in\balexr$.
Let $a\in\R$ and prove $\Phi$ is left continuous  at $a$.
Suppose $F(a)>G(a)$.
Given $\epsilon>0$ there  is  $\delta>0$ such that
$|F(x)-F(a)|<\epsilon$, $|G(x)-G(a)|<\epsilon$
and $F(x)>G(x)$
whenever $x\in(a-\delta,a)$.
For such $x$, $|\Phi(x)-\Phi(a)|=|F(x)-F(a)|<\epsilon$.
If $F(a)=G(a)$ then 
$|\Phi(x)-\Phi(a)|\leq \max(|F(x)-F(a)|,|G(x)-G(a)|)<\epsilon$.
Therefore, $\Phi$ is left continuous on $(-\infty,\infty]$.
For $x\in(-\infty,1/\delta)$ we can assume $\max(|F(x)|, |G(x)|)<\epsilon$.
Therefore, $|\Phi(x)|<\epsilon$.  Similarly, $\Phi$ has a right
limit at each point so that $\Phi\in\balexr$.  
Similarly with the infimum.  Hence,
$\Phi,\Psi\in\balexr$.

The following properties follow immediately from  the definition.
If $F\leq G$ then for all $H\in\balexr$ we have $F+H\leq G+H$.
If $F\leq G$ and $a\geq 0$ then $aF\leq aG$.  If $|F|\leq |G|$
then $\|F\|_\infty\leq \|G\|_\infty$.  Hence, $\balexr$ is a
Banach lattice.

(b), (c) First we show that
$\alexr$ is closed under the operations $f\vee g$ and
$f\wedge g$.  For $f,g\in\alexr$, we have
$f\vee g  =  \sup(f,g)$. This is $h$ such that
$f\preceq h$,  $g \preceq h$, and if $f\preceq\htilde$, $g\preceq\htilde$, then 
$h\preceq\htilde$.
This last statement is equivalent to 
$F\leq H$,  $G \leq H$, and if $F\leq\Htilde$, $G\leq\Htilde$, then 
$H\leq\Htilde$.
But then $H=\max(F,G)$ and $h=H'$ so $f\vee g=(F\vee G)'\in\alexr$.
Similarly, $f\wedge g=(F\wedge G)'\in\alexr$.
And, $|F'|=F'\vee(-F')=F'\vee(-F)'=(F\vee(-F))'=|F|'$.  The proofs that $(F^+)'=(F')^+$ and $(F^-)'=(F')^-$ are similar.
 
If $f,g\in\alexr$ and $f\preceq g$ then $F\leq G$.  Let $h\in\alexr$.  Then, $F+H\leq G+H$.
But then $(F+H)'=F'+H'=f+h\preceq g+h$.  If  $k\in\R$ and $k\geq 0$
then $(kF)'=kF'=kf$ so $kf\preceq kg$.  And, if $|f|\preceq |g|$ then 
$|F|'\preceq |G|'$ so $|F|\leq |G|$, i.e., $|F(x)|\leq |G(x)|$ for all $x\in\R$.
Then $\|f\|=\|F\|_\infty\leq \|G\|_\infty=\|g\|$. And, $\alexr$ is a
Banach lattice that is isomorphic to $\balexr$.

(d) Note that $|\int_{(a,b)}f|=|F(b-)-F(a+)|\geq | |F(b-)|-|F(a+)||
=|\int_{(a,b)}|F|'|=|\int_{(a,b)}|F'||$.  Similarly for other intervals.
The other parts of
(d) and (e) follow from (c) and the definitions.
(f) The real-valued functions on any set form a distributed lattice
due to inheritance from $\leq$ in $\R$.  Therefore, $\balexr$ is a distributed
lattice and hence so is $\alexr$. See \cite[p.~484]{maclanebirkhoff} for an
elementary
proof and for another property of distributed lattices.  (g) Modularity
is also inherited from $\leq$ in $\R$ via $\balexr$.
(h) We have $F',G'\in\alexr$ with respective primitives $\Phi_F,\Phi_G\in
\balexr$ given by $\Phi_F(x)=F(x-)-F(-\infty)$ and
$\Phi_G(x)=G(x-)-G(-\infty)$.  The definition of order then gives 
\eqref{lattice1}.  The relations $F(x\pm)=\lim_{y\to x^\pm}F(y)=
\lim_{y\to x^\pm}F(y-)$ then give \eqref{lattice2}.
\qed\\

For the
function $f(t)=H_1(t)\sin(t^2)$ we have
$f^+=|f|=f$ and $f^-=0$.

Notice that the definition of order allows us to integrate both sides of
$f\preceq g$ in
$\alexr$ to get $F\leq G$ in $\balexr$.  
The isomorphism allows us to differentiate
both sides of 
$F\leq G$ in $\balexr$ to get $F'\preceq G'$ in $\alexr$. 
If $F$ and $G$ are regulated functions on $\R$ with real limits
at $\pm\infty$ then the inequality $|F'|\preceq G'$ lets us
prove that $|F(x\pm)|-|F(-\infty)|\leq G(x\pm)-G(-\infty)$
for all $x\in\R$.  This is then a type of mean value theorem.
See \cite{koliha} where the inequality $|F'(x)|\leq G'(x)$
yields $|F(b)-F(a)|\leq G(b)-G(a)$ under the assumption
that $F$ is continuous or absolutely continuous and the
first inequality holds except on a countable set or set
of measure zero.  In \cite{koliha} $G$ is required to be increasing.

A lattice is {\it complete} if every subset that is bounded
above has a supremum.  But $\balexr$ is not complete.
Let $F_n(x)=H_1(x-1/n)|\sin(\pi/x)|$ and let $S=\{F_n\mid n\in\N\}$.
Then an upper bound for $S$ is $H_1$ but $\sup(S)(x)=H_1(x)|\sin(\pi/x)|$,
which is not regulated.
Hence, $\alexr$ is also not complete.

In this  section we have considered only the most elementary
lattice properties.  Other questions, such as the relation of
$\alexr$ and $\balexr$ to abstract $L$ spaces and abstract $M$
spaces, will be dealt with elsewhere.

\section{Topology and measure}\label{measure}
In this section we define a topology on $\R$ 
so that regulated functions
are continuous.  We then describe $\alexr$ in terms of finitely additive
measures.

The topology of half-open intervals or Sorgenfrey topology on
the real line is defined by taking a base to be the collection
of all intervals $(a,b]$ for all $-\infty<a<b<\infty$.  See,
for example, \cite[p.~156]{bauer}.  Call
the resulting topology $\tau_L$. Then $(\R, \tau_L)$ is
separable and  first countable but not second countable.
This topology is
finer than the usual topology on $\R$, hence it is a 
Hausdorff space. However,
it is not locally compact.  Each interval $(a,b]$ is also 
closed.  
Observe that
$[0,1]\subset (-1,0]\cup \cup_{n=1}^\infty (1/(n+1),1/n]$,
so $[0,1]$ is not compact in $\tau_L$.  In fact, each 
compact set is countable.

All functions in $\balexr$ are continuous in $(\R,\tau_L)$.
This follows from the fact that every regulated function is the uniform
limit of a sequence of step functions \cite[\S7.13]{mcleod}.  Hence it is only necessary
to consider $H_1$.  But $H_1^{-1}((0.5, 1.5))=(0,\infty)\in\tau_L$
and $H_1^{-1}((-0.5, 0.5))=(-\infty,0]\in\tau_L$.
Notice that right continuous functions need not be continuous in 
$(\R, \tau_L)$. For example, $H_2^{-1}((0.5,1.5))=[0,\infty)\not\in\tau_L$.
See Example~\ref{examples}(e) for the definition of
$H_2$.  Functions such as $f(x)=\sin(1/x)$ for $x>0$ and $f(x)=0$, otherwise,
and $g(x)=1/x$ for $x>0$ with $g(x)=0$, otherwise, are continuous
in $(\R,\tau_L)$, i.e., left continuous functions are continuous as
functions from $(\R, \tau_L)$ to $\R$ with the usual topology.

If $X$ is a nonempty set then an {\it algebra} on  
$X$ is a collection of
sets $\algebra\subset {\cal P}(X)$ such that (i) $\varnothing,X\in \algebra$
and if $E,F\in\algebra$ then (ii) $E\cup F\in\algebra$ and (iii)
$E\setminus F\in\algebra$.  Since $E\setminus F=(X\setminus F)\setminus
(X\setminus E)$, (iii) can be replaced with $X\setminus E\in\algebra$.
By de Morgan's laws, $\algebra$ is also closed under intersections.
Hence, $\algebra$ is closed under finite unions and intersections.
A set $E\subset\R$ is a $\bv$ set if $\chi_E\in\bv$. 
If $\algebra$ is an algebra then $\nu\fn\algebra\to\R$ is a
{\it finitely additive measure} if whenever $E,F\in\algebra$ such that
$E\cap F=\varnothing$ then $\nu(E\cup F)=\nu(E)+\nu(F)$.
Notice that $\nu(\varnothing)=0$.
We have the following results.
\begin{prop}
{\rm (a)} The $\bv$ sets form an algebra over $\R$.
{\rm (b)} If $f\in\alexr$ define $\nu_f(\varnothing)=0$ and 
$\nu_f(E)=\intinf f\chi_E$ for a $\bv$ set
 $E$. Then $\nu_f$ is a finitely additive measure on 
$\bv$ sets.  If $S$ is a $\bv$ set then $|f(S)|\leq \|f\|(1+ V\chi_S)$.
\end{prop}

\bigskip
\noindent
{\bf Proof:}
(a) Note that
$\chi_\varnothing =0\in\bv$ and $\chi_\R=1\in\bv$.  If $E$ and $F$ are
$\bv$ sets
then $\chi_{\R\setminus E}=1-\chi_E\in\bv$.  And,
$\chi_{E\cup F}=\chi_{\R\setminus[(\R\setminus E)\cap(\R\setminus F)]}
=1-\chi_{(\R\setminus E)\cap(\R\setminus F)}=1-\chi_{\R\setminus E}
\chi_{\R\setminus F}=1-[1-\chi_E][1-\chi_F]=\chi_E+\chi_F-\chi_E\chi_F\in\bv$.  
(b) Since the
functions of bounded variation are multipliers for $\alexr$ we have
that $\nu_f$ is a real-valued function on $\bv$ sets.  
If $E$ and $F$ are disjoint $\bv$ sets then $\chi_{E\cup F}=\chi_E+\chi_F$
and
$\nu_f(E\cup F)=\intinf f\chi_{E\cup F}=\intinf f(\chi_E+\chi_F)=\intinf f\chi_E
+\intinf f\chi_F =\nu_f(E)+\nu_f(F)$.  We have
$|f(S)|=|\intinf f\chi_S|\leq \|f\|(1+V\chi_S)$ using the H\"older
inequality Theorem~\ref{holder}. \qed\\

A finitely additive measure $\nu$ on algebra $\algebra$
is {\it finite} if $\sup_{E\in\algebra}|\nu(E)|<\infty$.  As finitely
additive measures, elements of $\alexr$ need not be finite.  For example,
if $g(x)=\sin(x)/x$ for $x\not=0$ then $f:=T_g\in\alexr$.  Let
$E_n=\cup_{k=0}^n[2k\pi,(2k+1)\pi]$. Then $\nu_f(E_n)\to\infty$
as $n\to\infty$.  We have the following connection between absolute integrability
and the
finitely additive measures that are finite in $\alexr$.

\begin{prop}\label{finitelyaddivenbv}
Let $f\in\alexr$.  Then $\nu_f$ is finite if and only if $F\in\nbv$.
\end{prop}

\bigskip
\noindent
{\bf Proof:} Suppose $F\in\nbv$ and $E$ is a $\bv$ set.  Then
$|\nu_f(E)|  =  \left|\intinf\chi_E\,dF\right|\leq VF$.  Hence, 
$|\nu_f|\leq VF<\infty$.

Suppose
$|\nu_f|<\infty$.  Let $(x_i,y_i)$ be disjoint intervals.  Then
$\sum|F(x_i)-F(y_i)|=\nu_f(\cup[x_i,y_i))$.  Therefore, $F\in\nbv$.
\qed\\

The $\bv$ sets do not form a $\sigma$-algebra.  For example,
the set $\cup_{n=1}^\infty[2n,2n+1]$ is not a $\bv$ set.

Let $F\fn\Rbar\to\R$ be any function.  Let $I$ be an interval with endpoints
$-\infty\leq a<b\leq\infty$.  Define $\nu(\varnothing)=0$ and
$\nu(I)=F(b)-F(a)$.  Then $\nu$ is a finitely additive measure on $\bv$
sets.  But $F$ need not be regulated so $\alexr$ does not contain all
finitely additive measures on $\bv$ sets.

\section{Convergence theorems}
There are different modes of convergence in $\alexr$.  If $\{f_n\}\subset
\alexr$ then $f_n\to f\in\alexr$ {\it strongly} if $\|f_n-f\|\to 0$.
The convergence is  {\it weak in} $\D$ if $\langle f_n-f,\phi\rangle
=\int_{-\infty}^\infty (f_n-f)\phi\to 0$ for all $\phi\in\D$
and the convergence is  {\it weak in} $\bv$ if 
$\int_{-\infty}^\infty (f_n-f)g\to 0$ for all $g\in\bv$.  Clearly,
strong convergence implies weak convergence in $\bv$ (Theorem~\ref{eqnorms}),
which
implies weak convergence in $\D$.  We would like conditions
under which $\int_{-\infty}^\infty f_n\to \int_{-\infty}^\infty f$.
Certainly weak convergence in $\bv$ is sufficient, take $g=1$.  Weak 
convergence in $\D$ is not sufficient.  For example, let $f_n=\tau_{n}\delta$,
for which $F_n(x)=H_1(x-n)$.  Then $\{f_n\}$ converges weakly in
$\D$ to $0$ but $F_n(\infty)=1$.  

Strong convergence is equivalent to
uniform convergence of the sequence of primitives.
\begin{theorem}\label{uniform}
Let $\{f_n\}\subset\alexr$ and let $\{F_n\}\subset\balexr$ be the
respective primitives.  Suppose $F\fn\Rbar\to\R$ and
$F_n\to F$ on $\Rbar$.  (a) $F_n\to F$ uniformly on $\Rbar$ if
and only if $\|f_n-f\|\to 0$.  (b) If $F_n\to F$ uniformly on $\Rbar$
then $F'\in\alexr$, $f_n\to F'$ strongly and $\intinf f_ng\to\intinf F'g$
for each $g\in\bv$.  In particular, $\int_If_n\to\int_IF'$ for each
interval $I\subset\R$.
\end{theorem}
Part (b) follows from the H\"older inequality.

If $\{F_n\}$ is  a sequence of continuous functions that converges
uniformly to function $F$  then $F$ is  continuous. 
A necessary and sufficient condition for $F$ to be continuous is
that the convergence be quasi-uniform.  Because of our compactification
of $\Rbar$ (Remark~\ref{compactification}), Arzel\`a's theorem
applies.  See 
\cite[pp.~268]{dunfordschwartz}.
We have a similar criteria
for regulated functions.

\begin{theorem}\label{qur}
Let each function $F_n\fn\Rbar\to\R$ be regulated on $\Rbar$.  
Suppose $F_n\to F$ at each point in $\Rbar$.   We require
$F(\pm\infty)\in\R$ and
$F_n(\pm\infty)\in\R$ for each $n\in\N$ but do not require
$F_n(\pm\infty)=\lim_{x\to\pm\infty}F_n(x)$.
Let $\epsilon>0$.
Suppose 
that for each $a\in\Rbar$ and each $N\in\N$ there exist $n\geq N$ and 
$\delta>0$ such that
if $x\in(a-\delta,a+\delta)$ then $|F_n(x)-F(x)|<\epsilon$,
for $a\in\R$.
For $a=-\infty$ we require $x\in[-\infty,-1/\delta)$.
For $a=\infty$ we require $x\in(1/\delta,\infty]$.
Then
$F$ is regulated on $\Rbar$.
\end{theorem}
\bigskip
\noindent
{\bf Proof:}
Let $\epsilon>0$.  Write 
$$
|F(x)-F(y)|\leq |F(x)-F_n(x)|+|F(y)-F_n(y)|+|F_n(x)-F_n(y)|.
$$
We have $n\geq 1$ and $\delta_n>0$ such that if $x\in(1/\delta_n,\infty]$
then $|F(x)-F_n(x)|<\epsilon$.  Each $F_n$ has a limit at $\infty$
so there is $\eta_n>0$ such that if $x,y\in(1/\eta_n,\infty)$ then
$|F_n(x)-F_n(y)|<\epsilon$.  Take $\delta=\min(\delta_n, \eta_n)$.
If $x,y\in(1/\delta,\infty)$ then $|F(x)-F(y)|<3\epsilon$.  Hence,
$\lim_{x\to\infty}F(x)$ exists.  Similarly, $F$ has a limit at $-\infty$.

The proof that $F$ has a left limit at $a\in\R$ is similar. Now
the intervals become
$x\in(a-\delta_n,a)$ and $x,y\in(a-\eta_n,a)$ and finally
$x,y\in(a-\delta,a)$.  Similarly for the right limit.\qed

See \cite[Proposition~3.6]{frankova} for a another sufficient 
condition on $\{F_n\}$
(bounded $\varepsilon$-variation)
that ensures $F$ is regulated.
\begin{corollary}\label{Brconv}
If $\{F_n\}\subset\balexr$ then $F\in\balexr$.
Now, as usual for functions in $\balexr$, we define 
$F_n(\pm\infty)=\lim_{x\to\pm\infty}F_n(x)$.
Let $f_n=F_n'$ and $f=F'$. Then for each $x\in(-\infty,\infty]$ we have
$\int_{(-\infty,x)}f_n\to\int_{(-\infty,x)}f$.
\end{corollary}
\bigskip
\noindent
{\bf Proof:}
Let $\epsilon>0$.  For $a\in\R$, write 
$$
|F(x)-F(a)|\leq |F(x)-F_n(x)|+|F_n(x)-F_n(a)|+|F_n(a)-F(a)|.
$$
Since $F_n(a)\to F(a)$ we have $N_a\in\N$ such that if $n\geq N_a$
then $|F_n(a)-F(a)|<\epsilon$.  We now have existence of 
$n\geq N_a$ and $\delta_n>0$
such that if $x\in(a-\delta_n,a]$ then $|F(x)-F_n(x)|<\epsilon$.
And, each $F_n$ is left continuous at $a$ so there is $\eta_n>0$ such
that if $x\in(a-\eta_n,a]$ then $|F_n(x)-F_n(a)|<\epsilon$.  Take
$\delta=\min(\delta_n,\eta_n)$.  If $x\in(a-\delta,a]$ then
$|F(x)-F(a)|<3\epsilon$ so $F$ is left continuous at $a$.
Similarly, $\lim_{x\to-\infty}F(x)=0=F(-\infty)$ and
$\lim_{x\to\infty}F(x)=F(-\infty)\in\R$.\qed

The Sorgenfrey topology of Section~\ref{measure} makes each
function in $\balexr$ continuous.  However, no
interval in $\Rbar$ is compact in this topology.  Hence,
Arzel\`a's theorem (\cite[pp.~268]{dunfordschwartz}), 
establishing that quasi-uniform convergence
is a necessary and sufficient condition for the limit of a 
sequence of continuous functions to be continuous, is not 
applicable.  We do not know necessary and sufficient conditions
under which a sequence in $\balexr$ will converge to   
a function in $\balexr$.  However, Theorem~\ref{uniform} and
Theorem~\ref{theoremconvergencebv} give
a sufficient condition while Theorem~\ref{qur} and Theorem~\ref{ang8} 
with their corollaries
give conditions under under which left continuity is preserved.

\begin{example}
{\rm
The example $f_n=\tau_{n}\delta$ in the first paragraph of this 
section shows the condition
at $\infty$ cannot be dropped.  For then we have
$F_n(x)=H_1(x-n)$.  For each $x\in\R$ we have $F_n(x)\to 0$ but
$F_n(\infty)=1$ so $F(x)=0$ for $x\in[-\infty,\infty)$ and $F(\infty)=1$.
Hence, $F\not\in\balexr$.  Although $f_n\to 0$ weakly in $\D$ we have
$\intinf f_n=F_n(\infty)=1\not\to 0$.  Note that if $n<x<\infty$ then
$|F_n(x)-F(x)|=1$ so the condition at infinity in Theorem~\ref{qur}
is not satisfied.
}
\end{example}

Weak convergence in $\D$ of $f_n$ to $f$ is not sufficient for
$\{F_n\}$ to converge to a function in $\balexr$.
The following theorem gives conditions in  addition to weak 
convergence in $\D$ so that $\int_{(-\infty,x)}f_n\to\int_{(-\infty,x)}f$.

\begin{theorem}\label{ang8}
Let $\{f_n\}\subset\alexr$ and let
$F\fn\R\to\R$ be regulated and left continuous 
on $\R$ with real limits at $\pm\infty$.
Suppose $\{F_n\}$ is uniformly bounded
on each compact interval in $\R$ and $F_n\to F$ on $\Rbar$.  Then
$f_n\to F'$ weakly in $\D$ and $\int_{(-\infty,x)} f_n\to\int_{(-\infty,x)} F'$
for each $x\in(-\infty,\infty]$.
\end{theorem}

The ordering introduced in Section~\ref{lattice} restores absolute
convergence to the integral.  Using this order we can rephrase
part of the above conditions in terms of dominated convergence.
\begin{corollary}\label{dominatedconvergence}
Let $\{f_n\}\subset\alexr$ and let
$F\fn\R\to\R$ be regulated and left continuous 
on $\R$ with real limits at $\pm\infty$.
Suppose there is $g\in\alexr$ such that $|f_n|\preceq g$ for all
$n\geq 1$.
Suppose $f_n\to f$ weakly in $\D$ for some $f\in \Dp$.
Suppose $F_n\to F$ on $\Rbar$.
Then $f=F'\in\alexr$
and
$\int_{(-\infty,x)} f_n\to\int_{(-\infty,x)} f$
for each $x\in (-\infty,\infty]$.
\end{corollary}

The proofs are easy modifications of Theorems~8 and 9 in \cite{ang}.
See also Theorem~17 in \cite{talviladenjoy}.

\begin{example}
{\rm
Let $f_n=n\chi_{(0,1/n)}-\tau_{1/n}\delta$.
Then $F_n(x)=nx$ for $0\leq x\leq 1/n$ and $F_n(x)=0$, otherwise.
We have $F_n\to 0$ on $\Rbar$.  The convergence is not uniform, since
$F_n(1/n)=1$.  Theorem~\ref{uniform} is not applicable.
The convergence $F_n\to 0$ on $\Rbar$ satisfies the conditions of 
Corollary~\ref{Brconv}.
This then gives $\int_{(-\infty,x)}f_n\to 0$ for each $x\in (-\infty,\infty]$.
Note that $|F_n(x)|\leq 1$ so $\{F_n\}$ is uniformly bounded.
Theorem~\ref{ang8} then gives the same conclusion.
Note that $0\leq F_n\leq H_1$ so $|f_n|\preceq H_1'=\delta\in\alexr$.
For $\phi\in\D$ we have 
$\intinf f_n\phi=n\int_0^{1/n}\phi(x)\,dx - \phi(1/n)\to 0$ by continuity.
Hence, 
$f_n\to 0$ weakly and
Corollary~\ref{dominatedconvergence} also gives the same conclusion.
}
\end{example}

The following theorem  follows from the H\"older inequality.
See also \cite{tvrdy1989}, \cite{tvrdy2002}.
\begin{theorem}[Uniform bounded variation]\label{theoremconvergencebv}
Suppose $\{f_n\}\subset\alexr$, $f\in\alexr$, $\{g_n\}\subset\bv$ and
$g\in\bv$ such that $\|f_n-f\|\to 0$, $V(g_n-g)\to 0$ and
$g_n(a)\to g(a)$ for some $a\in\R$.  Then $\|g_n-g\|_\infty\to 0$
and $\intinf f_ng_n\to\intinf fg$.
\end{theorem}

\bigskip
\noindent
{\bf Proof:}
First note that 
\begin{eqnarray*}
|g_n(x)-g(x)| & \leq & |g_n(a)-g(a)| + |[g_n(x)-g(x)] -
[g_n(a)-g(a)]|\\
 & \leq & |g_n(a)-g(a)| + V(g_n-g)\\
 & \to & 0\quad \text{ as } n\to\infty.
\end{eqnarray*}
Hence, $\|g_n-g\|_\infty\to 0$.  Now use the H\"older inequality
(Theorem~\ref{holder})
to write
\begin{eqnarray*}
\left|\intinf f_ng_n-\intinf fg\right| & = & 
\left|\intinf f_n(g_n-g) + (f_n-f)g\right|\\
 & \leq & \|f_n\|\,\|g_n-g\|_{\bv}
+ \|f_n-f\|\,\|g\|_{\bv}\\
 & \to & 0\quad \text{ as } n\to\infty.\qed
\end{eqnarray*}
\begin{corollary}
Suppose $\{f_n\}\subset\alexr$ and $f\in\alexr$
such that $\|f_n-f\|\to 0$.
If $g\in\bv$ then
$\intinf f_ng\to\intinf fg$.
\end{corollary}
\begin{corollary}\label{corconvergencebv}
Suppose $f\in\alexr$, $\{g_n\}\subset\bv$ and
$g\in\bv$ such that $V(g_n-g)\to 0$ and
$g_n(a)\to g(a)$ for some $a\in\R$.  Then
$\intinf fg_n\to\intinf fg$.
\end{corollary}
Note that in Theorem~\ref{theoremconvergencebv} and the
two corollaries we also get convergence on each subinterval
of $\R$.

\begin{example}[Convolution]
{\rm Suppose $f\in\alexr$ and $g\in AC$ such that both limits
$\lim_{x\to\pm\infty}g(x)$ exist in $\R$.  The convolution
$(f\ast g)(x)=\intinf f(x-y)g(y)\,dy$ exists for all
$x\in\R$.  By the change of variables theorem, $f\ast g=g\ast f$.
(Use $G(y)=x-y$ in Theorem~\ref{theoremchangevariables}.)  By
the H\"older inequality, $\|f\ast g\|_\infty\leq
\|f\|(\|g\|_\infty + Vg)$.  For each $x\in\R$,
we have $g(x-\cdot)\in\bv$.  For each $y, z\in\R$, we have
$\lim_{x\to z}g(x-y)=g(z-y)$.  And, 
\begin{eqnarray*}
V\left(g(z-\cdot)-g(x-\cdot)\right) & = &
\|g'(z-\cdot)- g'(x-\cdot)\|_1\\
 & \to & 0\quad\text{as }z\to x \text{ by continuity in the } 
L^1 \text{ norm}.
\end{eqnarray*}
By Corollary~\ref{corconvergencebv},
$f\ast g$ is 
uniformly
continuous on $\R$.
}
\end{example}

\section{Acknowledgments}
\v{S}tefan Schwabik provided several references.


\begin{thebibliography}{99}
\bibitem{alexiewicz}
A. Alexiewicz, {\it Linear functionals on Denjoy--integrable functions},
Colloquium Math. {\bf 1}(1948), 289--293.
\bibitem{ang}
D.D. Ang, K. Schmitt and L.K. Vy,
{\it A multidimensional analogue of the Denjoy--Perron--Henstock--Kurzweil
integral}, 
Bull. Belg. Math. Soc. Simon Stevin {\bf 4}(1997), 355--371.
\bibitem{aye}
K.K. Aye and P.-Y. Lee, 
{\it The dual of the space of functions of bounded variation},
Math. Bohem. {\bf 131}(2006), 1--9.
\bibitem{bauer}
H. Bauer, {\it Measure and integration theory}, Berlin,
Walter de Gruyter, 2001.
\bibitem{bongiorno}
B. Bongiorno, L. Di Piazza and D. Preiss, {\it
A constructive minimal integral which includes Lebesgue integrable 
functions and derivatives}, J.~London Math. Soc. (2) {\bf 62}(2000), 117--126.
\bibitem{celidze}
V.G. \v{C}elidze and A.G. D\v{z}var\v{s}e\v{\i}\v{s}vili,
{\it The theory of the Denjoy integral and some applications}
(trans. P.S. Bullen),
Singapore, World Scientific, 1989.
\bibitem{dales}
H.G. Dales et al, {\it Introduction to Banach algebras, operators,
and harmonic analysis}, Cambridge, Cambridge University Press, 2003.
\bibitem{dunfordschwartz}
N. Dunford and J. Schwartz, {\it Linear
Operators, Part I}, New York, Interscience, 1957.
\bibitem{evans}
L.C. Evans and R.F. Gariepy, {\it Measure theory and fine properties
of functions}, Boca Raton, CRC Press, 1992.
\bibitem{frankova}
D. Fra\v nkova, {\it Regulated functions},  Math. Bohem.  {\bf 116}(1991), 
20--59.
\bibitem{friedlander}
F.G. Friedlander and M. Joshi,
{\it Introduction to the theory of distributions},
Cambridge, Cambridge University Press, 1999.
\bibitem{hildebrandt1934}
T.H. Hildebrandt, {\it On bounded linear functional operations},
Trans. Amer. Math. Soc. {\bf 36}(1934), 868--875.
\bibitem{hildebrandt1938}
T.H. Hildebrandt, {\it Linear operations on functions of bounded
variation}, Bull.  Amer. Math. Soc.  {\bf 44}(1938), 75.
\bibitem{hildebrandt1966}
T.H. Hildebrandt, 
{\it Linear continuous functionals on the space $(BV)$ with weak topologies},
Proc. Amer. Math. Soc.  {\bf 17}(1966), 658--664.
\bibitem{honig}
C.S. H\"onig, {\it Volterra Stieltjes-integral equations},
Amsterdam, North-Holland, 1975.
\bibitem{jeffery}
R.L. Jeffrey, {\it The theory of functions of a real variable},
New York, Dover, 1985.
\bibitem{kaltenborn}
H.S. Kaltenborn, {\it Linear functional operations on functions having
discontinuities of the first kind},
Bull. Amer. Math. Soc. {\bf 40}(1934), 702--708.
\bibitem{kannan}
R. Kannan and C.K. Krueger, {\it Advanced analysis on the real line},
New York, Springer--Verlag, 1996.
\bibitem{koliha}
J.J. Koliha, {\it Mean, meaner, and the meanest mean value theorem},
Amer. Math. Monthly {\bf 116}(2009), 356--361.
\bibitem{maclanebirkhoff}
S. Mac Lane and G. Birkhoff, {\it Algebra}, New York, Macmillan, 1979.
\bibitem{mauldin}
R.D. Mauldin, 
{\it Some effects of set-theoretical assumptions in measure theory},
Adv. in Math. {\bf 27}(1978), 45--62.
\bibitem{mcleod}
R.M. McLeod, {\it The generalized Riemann integral}, Washington,
The Mathematical Association of America, 1980.
\bibitem{pmikusinski1}
P. Mikusi\'nksi and K. Ostaszewski, {\it Embedding Henstock integrable
functions into the space of Schwartz distributions},
Real Anal. Exchange {\bf 14}(1988-89), 24--29.
\bibitem{pfeffer}
W. Pfeffer,
{\it Derivation and integration}, Cambridge, Cambridge University Press, 2001.
\bibitem{raju}
C.K. Raju, {\it Distributional matter tensors in relativity},
Proceedings of the Fifth Marcel Grossmann Meeting on General Relativity, Part~A
(D.G. Blair and M.J. Buckingham, eds.),
Singapore, World Scientific, 1989, pp.~419--422.
\bibitem{saks}
S. Saks, {\it Theory of the integral} (trans. L.C. Young), Warsaw,  Monografie
Matematyczne, 1937.
\bibitem{schwabik1973}
\v{S}. Schwabik, {\it On the relation between Young's and Kurzweil's 
concept of Stieltjes integral},  \v Casopis P\v est. Mat.  {\bf 98}(1973), 
237--251.
\bibitem{schwabik1992}
\v{S}. Schwabik, {\it Linear operators in the space of regulated functions},
Math. Bohem.  {\bf 117}(1992), 79--92. 
\bibitem{schwartz}
L. Schwartz,
{\it Th\`eorie des distributions},
Paris, Hermann, 1966.
\bibitem{talvilafourier}
E. Talvila, {\it Henstock--Kurzweil Fourier transforms}, 
Illinois J. Math. {\bf 46}(2002), 1207--1226.
\bibitem{talviladenjoy}
E.  Talvila, {\it  The distributional Denjoy  integral}, Real Anal. Exchange
{\bf 33}(2008), 51--82.
\bibitem{tvrdy1989}
M. Tvrd\'y, {\it Regulated functions and the Perron--Stieltjes integral},
\v Casopis P\v est. Mat.  {\bf 114}(1989),  187--209. 
\bibitem{tvrdy1996}
M. Tvrd\'y, 
{\it Linear bounded functionals on the space of regular regulated functions},
Tatra Mt. Math. Publ.  {\bf 8}(1996), 203--210. 
\bibitem{tvrdy2002}
M. Tvrd\'y,  {\it Differential and integral equations in the space of 
regulated functions},  Mem. Differential Equations Math. Phys.  {\bf 25}(2002), 1--104. 
\bibitem{vladimirov}
V.S. Vladimirov, {\it Methods of the theory of generalized functions}, London,
Taylor and Francis, 2002.
\bibitem{ziemer}
W.P. Ziemer, {\it Weakly differentiable functions}, New York, Springer--Verlag,
1989.
\end{thebibliography}
\end{document}